\documentclass[journal]{IEEEtran}

\usepackage{graphicx}
\usepackage{caption}

\usepackage{algorithm}
\usepackage[noend]{algpseudocode}

\usepackage{multicol}
\usepackage{amsmath,amssymb}
\usepackage{amsfonts}
\usepackage{textcomp}
\usepackage{psfrag}
\usepackage{multimedia}
\usepackage{enumitem}

\newtheorem{theorem}{Theorem}[section]
\newtheorem{Proposition}[theorem]{Proposition}

\newtheorem{Corollary}[theorem]{Corollary}
\newtheorem{Lemma}[theorem]{Lemma}
\newtheorem{Definition}[theorem]{Definition}
\newtheorem{Remark}[theorem]{Remark}
\newtheorem{Example}[theorem]{Example}
\newtheorem{Assumption}[theorem]{Assumption}

\usepackage[usenames,dvipsnames]{color}
\definecolor{wheat}{rgb}{0.96,0.87,0.70}
\definecolor{mario}{rgb}{0.8,0.8,1}
\definecolor{seb}{rgb}{0.8,1,0.8}
\definecolor{myGreen}{rgb}{0.,0.8,0.0}

\definecolor{darkgreen}{rgb}{0,0.6,0}

\def\R{\mathbb{R}}
\def\eps{\epsilon}

\newcommand {\matr}[2]{\left[\begin{array}{#1}#2\end{array}\right]}

\newcounter{lastnote}

\usepackage{tikz}
\usetikzlibrary{calc,arrows,positioning}

\begin{document} 

\title{Stabilization of Strictly Pre-Dissipative Receding Horizon Linear Quadratic Control by Terminal Costs}

\author
{Mario Zanon, Lars Gr\"une
	\thanks{Mario Zanon is with the IMT School for Advanced Studies Lucca, Italy. }
	\thanks{Lars  Gr\"une is with the University of Bayreuth, Germany.}
	\thanks{}
}

\IEEEtitleabstractindextext{
	\begin{abstract}
		Asymptotic stability in receding horizon control is obtained under a strict pre-dissipativity assumption, in the presence of suitable state constraints. In this paper we analyze how terminal constraints can be replaced by suitable terminal costs. We restrict to the linear-quadratic setting as that allows us to obtain stronger results, while we analyze the full nonlinear case in a separate contribution.
	\end{abstract}
	
}

\maketitle

\IEEEdisplaynontitleabstractindextext

\IEEEpeerreviewmaketitle

\section{Introduction}

Receding horizon control (often used synonymously with model predictive control) is a control technique in which a finite horizon optimal control problem is solved in each time step and only the first element of the resulting optimal control sequence is used, while in the next time step the state is measured and the problem solved again to update the control~\cite{Rawlings2017,Gruene2017}. Under suitable stabilizability and regularity conditions, this scheme yields a practically asymptotically stable closed loop if the system is strictly dissipative with supply function defined via the stage cost of the finite horizon optimal control problem \cite[Chapter 8]{Gruene2017}. In this case, we call the optimal control problem {\em strictly dissipative}. Here, the size of the ``practical'' neighborhood of the equilibrium to which the closed-loop solution converges is determined by the length of the finite optimization horizon. True (as opposed to practical) asymptotic stability can be achieved by using suitable terminal constraints and costs, see \cite{DiAR10,Amrit2011a} or Theorem 8.13 in \cite{Gruene2017}. In these approaches the terminal cost is typically a local control Lyapunov function for the system and the terminal constraints are needed because the design of a global control Lyapunov function is usually a very difficult task. As a simpler alternative, it was shown in \cite{Zanon2018a} that linear terminal costs can also be used to obtain true asymptotic stability.

The strict dissipativity property that is at the heart of all these results requires the existence of a so-called storage function $\lambda$ mapping the state space into the reals. It is a strengthened version of the system theoretic dissipativity property introduced by Willems in his seminal papers \cite{Willems1972a,Willems1972b} and also featured in his slightly earlier paper \cite{Willems1971} on linear quadratic optimal control and the algebraic Riccati equation. Readers familiar with Lyapunov's stability theory can see the storage function $\lambda$ as a generalization of a Lyapunov function. However, unlike Lyapunov functions, $\lambda$ need not attain nonnegative values. However, it must be bounded from below, and this property is crucial for deriving the (practical) stability properties for receding horizon control cited above.

For generalized linear quadratic problems (by which we mean problems with linear dynamics and a cost function containing quadratic and linear terms) with state space $\R^{n_x}$, a standard construction for a storage function results in a function of the form $\lambda(x)= x^TPx + \nu^Tx$, for $P\in\R^{n_x\times n_x}$ and $\nu\in\R^{n_x}$, see \cite[Proposition 4.5]{DGSW14}. Clearly, such a function $\lambda$ is in general {\em not} bounded from below and Example~2.3 in \cite{DGSW14}, which we also present as Example \ref{ex:simple_lqr2}, below, shows that storage functions unbounded from below may occur even for very simple scalar problems. While the potential unboundedness of $\lambda$ has been handled somewhat informally in~\cite{DGSW14}, later in~\cite{Gruene2018} the variant of strict dissipativity with storage function not bounded from below has been termed {\em strict pre-dissipativity}. For strictly pre-dissipative problems, one way to restore strict dissipativity and thus (practical) asymptotic stability is to suitably restrict the state space by means of state constraints, e.g., to a compact set, on which $\lambda$ is bounded from below. 

Since such a restriction of the state space may not always be desirable, in this paper we will look at an alternative way to regain (practical) asymptotic stability. More precisely, we want to answer the following question: Given a receding horizon control scheme with strictly pre-dissipative optimal control problem, can we add a simple terminal cost that guarantees (practical) asymptotic stability? Here ``simple'' means that we don't want to design control Lyapunov function terminal costs but terms that are easier to compute. Thus, the answer we are looking for is conceptually similar to the one given in \cite{Zanon2018a}, but instead of the linear terminal cost from this reference we will see that a quadratic terminal cost is the right object here. We emphasize that no terminal constraints are needed in this approach. While we consider the general nonlinear case in the companion paper~\cite{Gruene2024}, in this paper we focus on the linear-quadratic case, for which we are able to provide stronger results, which we also connect to the generic results for the nonlinear case.

The analysis in this paper will deal with linear quadratic problems. After delivering the main results, we will briefly comment on the connections with the general nonlinear case. For this linear-quadratic problems, the question about asymptotic stability is closely linked with the existence of particular solutions to algebraic Riccati equations. Hence, we will make ample use of results from this area. The remainder of this paper is organized as follows. In Section~\ref{sec:statement} we introduce the problem and provide some definitions. We provide preliminary results in Section~\ref{sec:preliminary}, which are instrumental for our main results, delivered in Section~\ref{sec:main}. We discuss the connections with other results available in the literature in Section~\ref{sec:nonlinear}. We illustrate our theory with a numerical example in section~\ref{sec:simulations} and we draw our conclusions in Section~\ref{sec:conclusions}.

\section{Problem Statement}
\label{sec:statement}

We consider discrete-time systems of the form
\begin{align}
	\label{eq:system}
	x_{k+1} = A x_k + B u_k,
\end{align}
where $x\in\mathbb{R}^{n_x}$ and $u\in\mathbb{R}^{n_u}$ denote the states and controls respectively. 

Model predictive control consists in minimizing a given stage cost $\ell:\R^{n_x}\times\R^{n_u}\to\R$ over a given finite prediction horizon $N$, possibly subject to constraints and with the addition of a terminal cost. The receding horizon optimal control problem (RH-OCP) reads
\begin{subequations}
	\label{eq:empc}
	\begin{align}
		\min_{x_0,u_0,\ldots,x_N} \ \ & \sum_{k=0}^{N-1} \ell(x_k,u_k) + V^\mathrm{f}(x_N) \label{eq:empc_cost}\\
		\mathrm{s.t.} \ \ & x_0 = \hat x_j, \label{eq:empc_ic}\\
		&x_{k+1} = A x_k + B u_k,  && k \in\mathbb{I}_0^{N-1}, \label{eq:empc_dyn}\\
		&C x_k + D u_k \leq 0,  && k \in\mathbb{I}_0^{N-1}, \label{eq:empc_pc}
	\end{align}
\end{subequations}
where $h:\R^{n_x}\times\R^{n_u}\to\R^l$ defines the state and input constraints and inequality \eqref{eq:empc_pc} is to be understood componentwise. 

While we study the nonlinear case in the companion paper~\cite{Gruene2024}, in this paper we focus on the linear-quadratic case, i.e., 
\begin{subequations}
	\label{eq:lq_system}
	\begin{align}
		\ell(x,u) &= \matr{c}{x\\u}^\top H \matr{c}{x\\u}, \label{eq:lq_system:cost} \\
		V^f(x) &= x^\top P^\mathrm{f} x. \label{eq:lq_system:terminal_cost} 
	\end{align}
\end{subequations}
Note that, as the case of additional constant terms in~\eqref{eq:system} and additional linear terms in~\eqref{eq:lq_system:cost}-\eqref{eq:lq_system:terminal_cost} has been fully analyzed in~\cite{Zanon2018a}, in this section we assume for simplicity that these quantities are zero. 

Starting in the initial value $\hat x_0$ at time instant $j=0$, at every time instant $j\ge 0$ the state $\hat x_j$ is measured, Problem~\eqref{eq:empc} is solved, and the first optimal input $u_0^\star$ is applied to the system to obtain
\begin{align}
	\label{eq:clsystem}
	\hat x_{j+1} = A\hat x_j + B u_0^\star. 
\end{align}%
This procedure is repeated iteratively for all $j\ge 0$.

In this paper, we are interested in obtaining stability properties of the closed loop system~\eqref{eq:clsystem}. While stability results are abundant for the case of suitably formulated terminal constraints and Lyapunov function terminal costs~\cite{Diehl2011, Angeli2012a,Faulwasser2018a,Muller2015,Muller2013a}, we focus next on the case of no terminal constraint. This case has been analyzed, e.g., in~\cite{Faulwasser2018,Zanon2018a,Grune2013a,Muller2016}, where we can further distinguish between formulations without terminal cost and formulations with simple terminal costs that need not be Lyapunov functions, which are usually difficult to design. This last approach has in particular been taken in~\cite{Faulwasser2018,Zanon2018a} by using a linear terminal cost and the present paper can be seen as a continuation of this research. As in these references, our analysis is based on dissipativity concepts. 

We define the infinite-horizon Optimal Control Problem (OCP) related to~\eqref{eq:lq_system}, known as the Linear Quadratic Regulator (LQR), as
\begin{subequations}
	\label{eq:lqr_terminal_cost}
	\begin{align}
		V(x_0)=\min_{u_0,u_1,\ldots} \ & \lim_{N\to\infty}\sum_{k=0}^{N-1} \ell(x_k,u_k) + x_N^\top P^\mathrm{f} x_N\\
		\mathrm{s.t.} \ \ & x_{k+1} = A x_k + B u_k.
	\end{align}
\end{subequations}
Under suitable assumptions, the LQR solution is characterized by a solution $P$ with associated $K$ of the Discrete Algebraic Riccati Equation (DARE)
\begin{subequations}
	\label{eq:DARE}	
	\begin{align}
		P &= Q + A^\top P A - (S^\top + A^\top P B)K, \\
		K &= (R+B^\top P B)^{-1} (S+B^\top P A), 
	\end{align}
\end{subequations}
which provides both the optimal feedback law
\begin{align*}
	F(x) = -Kx,
\end{align*}
defining the optimal control via $u_k^\star = -Kx_k^\star$, and the associated quadratic value function
\begin{align*}
	V(x) = x^\top P x. 
\end{align*}
For a finite horizon $N$, the cost-to-go and the optimal feedback law are  in this case time-varying and read respectively $V_k(x) = x^\top P_{N-k} x$, and $F_k(x) = -K_{N-k} x$, where
\begin{subequations}
	\label{eq:RE}	
	\begin{align}
		P_{n+1} &= Q + A^\top P_n A - (S^\top + A^\top P_n B)K_{n+1},  \label{eq:RE:P}\\
		K_{n+1} &= (R+B^\top P_n B)^{-1} (S+B^\top P_n A),  \label{eq:RE:K}
	\end{align}
\end{subequations}
defines $P_1,\ldots,P_N$ and $K_0,\ldots,K_N$ inductively with $P_0=P^\mathrm{f}$. Note that for an RH-OCP with data of the linear quadratic form \eqref{eq:lq_system} the receding horizon control in \eqref{eq:clsystem} is given by 
\[ u_0^\star = -K_N(\hat x_j).\]


In the context of this paper, it is paramount to clarify the fact that matrix $P$ defining the optimal value function of the infinite-horizon LQR problem and the (symmetric) solutions $P$ of the DARE do not necessarily coincide, as we explain next. 
It is well-known that the DARE has in general infinitely many symmetric solutions~\cite{Willems1971,Molinari1975a,Lancaster1986,Ran1993a}, which can be interpreted as solutions of different OCPs. Note that the DARE might also have non-symmetric solutions, but these are not relevant in the context of this paper. Therefore, whenever referring to a solution of the DARE we will implicitly restrict to the set of symmetric solutions.

We call a solution $P_\mathrm{s}$ with associated $K_\mathrm{s}$ {\em stabilizing}, if all eigenvalues $\mu$ of $A-BK_\mathrm{s}$ satisfy $|\mu|<1$, i.e., if $A-BK_{\mathrm{s}}$ is Schur stable. 
If there exists a stabilizing solution to the DARE, then it is unique, i.e., all other solutions are not stabilizing~\cite{Lancaster1986,Ran1993a}. 
Unfortunately, even in case the DARE does have a stabilizing solution, this one need not be the solution to the LQR, as proven by the next example.

\begin{Example}
	\label{ex:simple_lqr2}
	Consider an LQR problem defined by $A=2$, $B=1$, $Q=0$, $S=0$, $R=1$, and $P^\mathrm{f}=0$. 
	The DARE reads 
	\begin{align*}
		P = 4P - \frac{4P^2}{1+P}, && K= \frac{2P}{1+P},
	\end{align*}
	and has the two solutions
	\begin{align*}
		P = \{0,3\}, && K=\{0,1.5\}.
	\end{align*}
	Both solutions of the DARE correspond to the solution of an OCP. The first one is
	\begin{align*}
		\min_{u_0,u_1,\ldots} \ & \sum_{k=0}^\infty u_k^2 \\
		\mathrm{s.t.} \ \ & x_{k+1} = A x_k + B u_k, && k\in\mathbb{I}_0^\infty,
	\end{align*}
	which has the trivial solution $u_k^\star=0$, for all $k$ and corresponds to $P=0$, $K=0$, which does not stabilize the system. 
	The second OCP is
	\begin{align*}
		\min_{u_0,u_1,\ldots} \ & \sum_{k=0}^\infty u_k^2 \\
		\mathrm{s.t.} \ \ & x_{k+1} = A x_k + B u_k,  && k\in\mathbb{I}_0^\infty,\\
		& \lim_{k\to\infty} x_k =0,
	\end{align*}
	which corresponds to the stabilizing solution $P=3$ and $K=1.5$. Thus, the asymptotic state constraint $\lim_{k\to\infty} x_k =0$ leads to a stabilizing optimal solution.
	
	Note that for this system, the storage function $\lambda(x) = - x^2$ yields the rotated cost matrix
	\begin{align*}
		\bar H = \matr{ll}{3 & 2 \\ 2 & 2} \succ0,
	\end{align*}
	whose positive definiteness confirms that strict pre-dissipativity holds. Yet, asymptotic stability cannot be achieved unless suitable constraints are introduced. As already mentioned, we will prove in this paper that the constraints can be replaced by a suitably defined quadratic terminal cost. In the case of this example, one can see that any terminal cost $V^f(x) = x^\top P^\mathrm{f}x$ with $P^\mathrm{f}>0$ yields asymptotic stability. $\hfill\square$
\end{Example}

\section{Preliminary Results}
\label{sec:preliminary}

In this section we provide some useful results that we will exploit in order to prove our main contributions.

\subsection{Cost Rotation and Pre-Stabilization}
In order to derive our results, let us provide a few useful facts next.

\paragraph{Strict $(x,u)$-Pre-Dissipativity and Rotated Cost}
For a given matrix $\Lambda$ and storage function $\lambda(x)=x^\top \Lambda x$, we define the rotated cost 
\begin{subequations}
	\label{eq:rotated_LQcost}
	\begin{align}
		L(x,u) := \matr{c}{x\\u}^\top H_\Lambda \matr{c}{x\\u}, 
	\end{align}
	with 
	\begin{align}\label{eq:Hlambda}
		H_\Lambda &:= \matr{ll}{Q + \Lambda - A^\top \Lambda A & S^\top - A^\top \Lambda B \\ S - B^\top \Lambda A & R - B^\top \Lambda B}. 
	\end{align}
\end{subequations}

\begin{Definition}[Quadratic Strict $(x,u)$-Pre-Dissipativity]
	Quadratic strict $(x,u)$-pre-dissipativity holds if there exists a scalar $\epsilon>0$ such that
	\begin{align}
		\label{eq:quadratic_strict_dissipativity}
		L(x,u) \geq \epsilon (\|x\|^2 + \|u\|^2), 
		&& \text{or, equivalently,} && H_\Lambda \succ 0.
	\end{align}
\end{Definition}
It has been proven in \cite[Lemma 4.1]{Gruene2018} that for linear quadratic problems the more generic concept of strict pre-dissipativity implies quadratic strict pre-dissipativity. In this paper we restrict to the slightly stronger case of quadratic strict $(x,u)$-pre-dissipativity, in order to guarantee the existence of a solution to the DARE. In practice, this excludes singular LQ problems.

In contrast to the more common strict dissipativity concept, strict pre-dissipativity, introduced under this name in \cite{Gruene2018}, does not require the storage function $\lambda$ to be bounded from below. This implies that one cannot use arguments as, e.g., in~\cite{Grune2013a,Gruene2014} in order to conclude (practical) stability properties of the closed loop \eqref{eq:clsystem}, and in fact stability may fail to hold, as we will show by means of the following example.

\begin{Example}\label{ex:simple_lqr1}
	Consider the LQR problem from Example~\ref{ex:simple_lqr2}. 
	One easily sees that for any initial condition $x_0$ and any horizon $N$ the optimal control sequence is $u_k^\star \equiv 0$, as this is the only control that produces $0$ cost, while all other control sequences produce positive costs. This implies that system \eqref{eq:clsystem} becomes 
	\[ \hat x_{j+1} = 2 \hat x_j, \]
	for which the origin is exponentially unstable. Yet, one checks that this problem is quadratically strictly $(x,u)$-pre-dissipative with storage function $\lambda(x) = - c x^2$ for each $c\in]0,3[$. 
	This shows that quadratic strict $(x,u)$-pre-dissipativity does not imply asymptotic stability of the optimal equilibrium.  $\hfill\square$
\end{Example}

As already mentioned in the introduction and as also seen in this example, storage functions that are not bounded from below appear naturally already for linear quadratic problems. In order to achieve closed-loop stability, often a compact state constraint set is imposed, as compactness implies boundedness of the storage function provided it is continuous (which is often the case). For Example \ref{ex:simple_lqr1}, it was shown in \cite[Example 2.3]{DGSW14} that this indeed renders the origin practically asymptotically stable for the closed loop. Yet, imposing compact state constraints just for the sake of achieving stability may not always be desirable.
As we will prove in this paper, stability can be alternatively achieved by a suitably defined quadratic terminal cost.

The next lemma formalizes the interest in rotated costs of the form~\eqref{eq:rotated_LQcost} in the context of this paper.
\begin{Lemma}
	\label{lem:rotated_equivalence}
	For any finite horizon $N$ as well as for the infinite horizon problem, an LQ problem with stage cost matrix $H$ and terminal cost matrix $P_0$ yields the same feedback law $u_k=-K_N x_k$ as the rotated LQ problem with stage cost matrix $H_\Lambda$ and terminal cost matrix $P_0+\Lambda$. Moreover, the matrices defining the optimal value functions of the two problems satisfy $P_{\lambda,N} = P_N + \Lambda$.
\end{Lemma}
\begin{IEEEproof}
	The proof follows from the observation that
	\begin{align*}
		\phantom{.}\hspace{3em}&\hspace{-3em}\sum_{k=0}^{N-1} \matr{c}{x_k\\u_k}^\top H_\Lambda \matr{c}{x_k\\u_k} + x_N^\top (P_0+\Lambda) x_N\\
		= \ & \sum_{k=0}^{N-1} \matr{c}{x_k\\u_k}^\top H \matr{c}{x_k\\u_k} + x_N^\top P_0 x_N + \sum_{k=0}^{N-1} x_k^\top \Lambda x_k \\
		& - \sum_{k=0}^{N-1} (Ax_k+Bu_k)^\top \Lambda (Ax_k+Bu_k) + x_N^\top \Lambda x_N \\
		= \ &\sum_{k=0}^{N-1} \matr{c}{x_k\\u_k}^\top H \matr{c}{x_k\\u_k} + x_N^\top P_0 x_N + \sum_{k=0}^{N-1} x_k^\top \Lambda x_k \\
		&- \sum_{k=0}^{N-1} x_{k+1}^\top \Lambda x_{k+1} + x_N^\top \Lambda x_N \\
		= \ &\sum_{k=0}^{N-1} \matr{c}{x_k\\u_k}^\top H \matr{c}{x_k\\u_k} + x_N^\top P_0 x_N + x_0^\top \Lambda x_0.
	\end{align*}
	Consequently, the costs of the two problems only differ by the constant term $x_0^\top \Lambda x_0$, hence their minimizers coincide. The same holds true if we consider the limit for $N\to\infty$ of the cost.
\end{IEEEproof}

This lemma entails that we can study the stability properties of any problem satisfying quadratic strict $(x,u)$-pre-dissipativity as an LQR problem with a positive-definite stage cost, provided that we suitably modify the terminal cost. However, while optimal control problems with positive-definite stage costs are always stabilizing (for sufficiently large $N$ and in the infinite-horizon case) if there is a positive-semidefinite terminal cost, this may no longer be true in the presence of a terminal cost which is not positive definite. Hence, since in the pre-dissipative case the rotated terminal cost matrix $P^\mathrm{f}+\Lambda$ can be indefinite even in case $P^\mathrm{f}\succ0$, quadratic strict $(x,u)$-pre-dissipativity does not necessarily entail that the LQR problem will be stabilizing. An example of such a situation is given once more by the system in Example~\ref{ex:simple_lqr2} with $P^\mathrm{f}=0$, for which the terminal cost matrix adapted to the rotated problem is $P^\mathrm{f} + \Lambda = -1$ and the corresponding optimal feedback $K=0$ is indeed not stabilizing.


\paragraph{Pre-Stabilized System}
Let us define $A_{\hat K}:=A-B{\hat K}$ and consider the following reformulation of problem~\eqref{eq:lqr_terminal_cost}, where we apply the change of variable $\bar u_k = u_k-{\hat K}x_k$:
\begin{subequations}
	\label{eq:lqr_terminal_cost2}
	\begin{align}
		\min_{\bar u_0,\bar u_1,\ldots} \ & \lim_{N\to\infty}\sum_{k=0}^{N-1} \ell(x_k,\bar u_k-{\hat K}x_k) + x_N^\top P^\mathrm{f} x_N\\
		\mathrm{s.t.} \ \ & x_{k+1} = A_{\hat K} x_k + B \bar u_k.
	\end{align}
\end{subequations}
In this case, we have that the stage cost matrix reads
\begin{subequations}
	\label{eq:pre_stabilized_cost}
	\begin{align}
		H_{\hat K} &= \matr{ll}{Q_{\hat K} & S_{\hat K}^\top \\ S_{\hat K} & R}, \\
		Q_{\hat K} &:= Q - S^\top {\hat K} - {\hat K}^\top S + {\hat K}^\top R {\hat K}, \\
		S_{\hat K} &:= S - R{\hat K}.
	\end{align}
\end{subequations}
With slight abuse of terminology, we will refer to this problem as pre-stabilized, even though $A_{\hat K}$ need not be Schur stable.

By Lemma~\ref{lem:rotated_equivalence} if the problem is strictly pre-dissipative then the optimal solution of~\eqref{eq:lqr_terminal_cost2} is also optimal for
\begin{subequations}
	\label{eq:lqr_terminal_cost2_rot}
	\begin{align}
		\hspace{-1em}\min_{\bar u_0,\bar u_1,\ldots} \ & \lim_{N\to\infty}\sum_{k=0}^{N-1} L(x_k,\bar u_k-{\hat K}x_k) + x_N^\top (P^\mathrm{f}+\Lambda) x_N\\
		\mathrm{s.t.} \ \ & x_{k+1} = A_{\hat K} x_k + B \bar u_k,
	\end{align}
\end{subequations}
where the rotated stage cost matrix reads
\begin{align*}
	H_{{\hat K},\Lambda} &:= \matr{ll}{Q_{\hat K} + \Lambda - A_{\hat K}^\top \Lambda A_{\hat K}& S_{\hat K}^\top - A_{\hat K}^\top \Lambda B \\ S_{\hat K} - B^\top \Lambda A_{\hat K} & R - B^\top \Lambda B} \succ 0.
\end{align*}
Finally, as we prove next, quadratic strict $(x,u)$-pre-dissipativity is not affected by pre-stabilizing the system.
\begin{Lemma}
	\label{lem:pre_stabilized_system}
	Assume that quadratic strict $(x,u)$-pre-dissipativity holds for a linear quadratic problem~\eqref{eq:lqr_terminal_cost}. Then for any $\hat K\in\R^{n_u\times n_x}$ quadratic strict $(x,u)$-pre-dissipativity also holds for the pre-stabilized problem~\eqref{eq:lqr_terminal_cost2}. Moreover, any solution of the DARE associated with the original problem~\eqref{eq:lqr_terminal_cost} is also a solution of the DARE associated with the pre-stabilized problem~\eqref{eq:lqr_terminal_cost2}.
\end{Lemma}
\begin{IEEEproof}
	With a few algebraic manipulations, one can show that
	\begin{align}
		M_{\hat K}^\top H M_{\hat K} &= H_{{\hat K}}, \\
		M_{\hat K}^\top G M_{\hat K} &= \matr{ll}{\Lambda - A_{\hat K}^\top \Lambda A_{\hat K} & - A_{\hat K}^\top \Lambda B \\ - B^\top \Lambda A_{\hat K} & - B^\top \Lambda B} =: G_{\hat K}, 
	\end{align}
	where 
	\begin{align*}
	G&:= \matr{ll}{\Lambda - A^\top \Lambda A & - A^\top \Lambda B \\ - B^\top \Lambda A & - B^\top \Lambda B}, &
	M_{\hat K} &:= \matr{ll}{I & 0 \\ -{\hat K} & I}.
	\end{align*}
	This entails $M_{\hat K}^\top H_\Lambda M_{\hat K}=H_{{\hat K},\Lambda}$. 
	Because $M_{\hat K}$ is full rank by construction, this yields that quadratic strict $(x,u)$-pre-dissipativity for the original problem, i.e., $H_\Lambda\succ0$, implies $H_{{\hat K},\Lambda}\succ0$, i.e., quadratic strict $(x,u)$-pre-dissipativity for the pre-stabilized problem.
	
	In order to prove the second claim, we write the DARE associated with the pre-stabilized problem~\eqref{eq:lqr_terminal_cost2}:
	\begin{align*}
		P &= Q_{\hat K} + A_{\hat K}^\top P A_{\hat K} - (S_{\hat K}^\top + A_{\hat K}^\top PB)K_{\hat K}, \\
		K_{\hat K}&= (R+B^\top PB)^{-1}(S_{\hat K} + B^\top PA_{\hat K}). 
	\end{align*}
	We then observe that
	\begin{align*}
		K_{\hat K}
		&=(R+B^\top PB)^{-1}(S + B^\top PA) - \hat K \\
		&= K-\hat K,
	\end{align*}
	such that, after few algebraic manipulations one obtains
	\begin{align*}
		P
		&= Q + A^\top P A - (S^\top + A^\top P B)K,
	\end{align*}
	i.e., the DARE associated with~\eqref{eq:lqr_terminal_cost}.
\end{IEEEproof}

We prove next some additional useful results.

\begin{Lemma}
	\label{lem:dare_pd}
	Assume that $(A,B)$ is controllable and quadratic strict $(x,u)$-pre-dissipativity holds. Then any symmetric solution $P$ of the DARE~\eqref{eq:DARE} satisfies $R+B^\top P B \succ 0$.
\end{Lemma}
\begin{IEEEproof}
	Because quadratic strict $(x,u)$-pre-dissipativity holds, by~\cite{Zanon2014d} a stabilizing solution exists, i.e., there exists a $P$ solving the DARE~\eqref{eq:DARE} with $R+B^\top P B \succ 0$.
	The proof then follows from~\cite[Theorem~2.5]{Lancaster1986}, (see also~\cite{Stoorvogel1998}), which states that, under the given assumption, if there exists one Hermitian solution $P$ such that $R+B^\top P B \succ 0$, then the condition holds for all Hermitian solutions. 
	%
\end{IEEEproof}	

This result can be then extended to the case of stabilizable systems as follows.

\begin{theorem}
	\label{thm:dare_pd}
	Assume that $(A,B)$ is stabilizable and quadratic strict $(x,u)$-pre-dissipativity holds. Then any symmetric solution $P$ of the DARE~\eqref{eq:DARE} satisfies 
	\begin{align}
		\label{eq:pd_dare_result}
		R+B^\top P B \succ 0.
	\end{align}
\end{theorem}
\begin{IEEEproof}
	Note that, since we assume stabilizability instead of controllability, the result of Lemma~\ref{lem:dare_pd} needs to be sharpened.

	Let us assume for simplicity and without loss of generality that the system is already in the canonical controllability form:
	\begin{align}
		A = \matr{ll}{A_{11} & A_{12} \\ 0 & A_{22}}, && B=\matr{ll}{B_1 \\ 0},
	\end{align}
	with $(A_{11},B_1)$ controllable. This allows us to separate the solution to the DARE in components. In particular, through few algebraic manipulations one can see that the DARE can be split in components, each defining one component of $P$, split as
	\begin{align*}
		P=\matr{ll}{P_{11} & P_{12} \\ P_{12}^\top & P_{22} }.
	\end{align*}
	Through few algebraic manipulations, one can see that the controllable part of the DARE reads
	\begin{align*}
		P_{11} &= Q_{11} + A_{11}^\top P_{11} A_{11} - (S_1^\top + A_{11}^\top P_{11} B_1)K_1, \\
		K_1 &= (R+B_1^\top P_{11} B_1)^{-1} (S_1+B_1^\top P_{11} A_{11}),
	\end{align*}
	i.e., it is a DARE in $P_{11}$ which is independent of $P_{12}$ and $P_{22}$.

	It is worth noting that the term $R+B^\top PB$ reads
	\begin{align*}
		R+\matr{ll}{B_{1}^\top & 0} \matr{ll}{P_{11} & P_{12} \\ P_{21} & P_{22}}\matr{ll}{B_{1} \\ 0 } 
		&= R+B_{1}^\top P_{11} B_{1},
	\end{align*}
	such that positive-definiteness of $R+B^\top PB$ can be assessed by considering the controllable part of the system only. 	
	The result then directly follows from Lemma~\ref{lem:dare_pd}.
\end{IEEEproof}

Since the computations in the proof only characterize part of the solution to the DARE, let us briefly discuss the full solution. By few algebraic manipulations one can see that: (a) the sub-equation defining $P_{12}$ does depend on $P_{11}$ but not on $P_{22}$, and the sub-equation is linear in $P_{12}$, such that the solution is unique; (b) once $P_{11}$ and $P_{12}$ have been solved for, also the sub-equation in $P_{22}$ is linear in $P_{22}$ and, consequently, has a unique solution.

Indeed, we have
\begin{align*}
	P_{12} &= \tilde Q_{12} + \tilde A_{12}^\top P_{12} A_{22},\\
	\tilde Q_{12} &:= Q_{12}  + A_{11}^\top P_{11} A_{12} - (S_2^\top + A_{11}^\top P_{11}B_{1}) K_{12} \\
	K_{12} &:= (R+B_{1}^\top P_{11} B_{1})^{-1}(S_1 + B_{1}^\top P_{11}A_{12} +  B_{1}^\top P_{12}A_{22}), \\
	\tilde A_{12}^\top &:= A_{12}^\top - (S_2^\top + A_{11}^\top P_{11}B_{1})(R+B_{1}^\top P_{11} B_{1})^{-1}B_{1}^\top.
\end{align*}
and, remembering that $P_{21}=P_{12}^\top$, 
\begin{align*}
	P_{22} &= \tilde Q_{22} +  A_{22}^\top P_{22}A_{22}, \\
	\tilde Q_{22} &:= Q_{22} + A_{12}^\top P_{11}A_{12} +  A_{22}^\top P_{21}A_{12} + A_{12}^\top P_{12} A_{22} \\
	&\hspace{7em}- (S_2^\top + A_{12}^\top P_{11}B_{1} +  A_{22}^\top P_{21}B_{1}) K_{22}, \\
	K_{22} &:= (R+B_{1}^\top P_{11} B_{1})^{-1} ( S_2 + B_{1}^\top P_{11}A_{12} +  B_{1}^\top P_{12}A_{22} ).
\end{align*}

	%

%
These results allow us to formulate the following conclusion.
\begin{Remark}
	The stability properties of any stabilizable linear system with quadratic stage cost satisfying the quadratic strict $(x,u)$-pre-dissipativity condition~\eqref{eq:quadratic_strict_dissipativity} can be studied by considering the corresponding fully controllable part with the corresponding rotated cost. Moreover, whenever matrix $A_{11}$ is singular, it will be convenient to pre-stabilize it so as to make it nonsingular. Since pre-stabilization does not alter the results, this approach remains general.
\end{Remark}


In the light of this remark, in the remainder of this section we make the following assumptions for simplicity and without loss of generality. 
\begin{Assumption}
	The linear system $(A,B)$ is controllable, $A$ is nonsingular and the stage cost matrix $H$ is positive definite, i.e., $H\succ0$.
\end{Assumption}

Before establishing our results, we need to provide some further definitions and preliminary results.

\subsection{The Reverse Discrete-Time Algebraic Riccati Equation}

As we will discuss next, among the symmetric solutions of the DARE, when they exist two of them play a fundamental role in the context of this paper: the stabilizing solution $P_\mathrm{s}$ with corresponding feedback $K_\mathrm{s}$ such that $A-BK_\mathrm{s}$ has all the eigenvalues strictly inside the unit circle; and the antistabilizing solution $P_\mathrm{a}$ with corresponding feedback $K_\mathrm{a}$ such that $A-BK_\mathrm{a}$ has all the eigenvalues strictly outside the unit circle. 

Unfortunately, even in case the DARE has a stabilizing solution, existence of the antistabilizing solution is not guaranteed in general. We will discuss next that the existence of $P_\mathrm{a}$ is guaranteed under the assumption that matrix $R-SA^{-1}B$ is nonsingular. In this case, $P_\mathrm{a}$ coincides with the stabilizing solution $\bar P_\mathrm{s}$ to the so-called   Reverse Discrete-time Algebraic Riccati Equation (RDARE). Since we will prove that $\bar P_\mathrm{s}$ exists whenever $P_\mathrm{s}$ exists, we will then also prove that it can be used as a substitute for $P_\mathrm{a}$ in case $P_\mathrm{a}$ does not exist.
%

The Reverse Discrete-time Algebraic Riccati Equation (RDARE) reads:
\begin{align*}
	\bar P &= \bar Q + \bar A^\top \bar P \bar A - (\bar S^\top + \bar A^\top \bar P \bar B) \bar K, \\
	\bar K &= (\bar R + \bar B^\top \bar P \bar B)^{-1}(\bar S + \bar B^\top \bar P \bar A),
\end{align*}
with 
\begin{align*}
	\bar A &:= A^{-1}, & \bar B &:= A^{-1} B, \\
	\bar Q &:= -\bar A^\top Q \bar A, & \bar S &:= S \bar A - \bar B^\top Q \bar A, \\
	\bar R &:= -R + S \bar B + \bar B^\top S^\top - \bar B^\top Q \bar B. \hspace{-10em}
\end{align*}
Note that we assumed without loss of generality that $A^{-1}$ exists.


\begin{Proposition}[{\cite[Proposition~2, Remark~2]{Ionescu1996}}]
	The DARE and RDARE share the same solutions if and only if $R-SA^{-1}B$ is nonsingular. In particular, if they exist, then $P_\mathrm{a}=\bar P_\mathrm{s}$ and $P_\mathrm{s}=\bar P_\mathrm{a}$, where $\bar P_\mathrm{a}$ is the antistabilizing solution of the RDARE. Moreover, if $R-SA^{-1}B$ is singular, then $R+B^\top \bar P B$ is singular for any symmetric $\bar P$ solving the RDARE.
\end{Proposition}
We observe that the second claim follows from the observation that 
\begin{equation*}
	R + B^\top \bar P B = (R-SA^{-1}B)^\top (\bar R + \bar B^\top \bar P \bar B)^{-1} (R-SA^{-1}B),
\end{equation*}
for any symmetric $\bar P$ solving the RDARE. This provides a clear intuition as to why the antisymmetric solution exists if and only if $R-SA^{-1}B$ is nonsingular.
Finally, note that also 
\begin{equation*}
	\bar R + \bar B^\top P \bar B = (R-SA^{-1}B)^\top (R + B^\top P B)^{-1} (R-SA^{-1}B)
\end{equation*}
holds. Consequently, if $R-SA^{-1}B$ is singular, then no solution $P$ of the DARE can be a solution $\bar P$ of the RDARE and viceversa.


\begin{Example}
	Consider the scalar system with $A=1$, $B=1$, and stage cost matrices $Q=1$, $R=1$, $S=1$. The DARE reduces to $P=1$, which corresponds to $K=1$, such that $A-BK=0$ is stable. Moreover, we have $\bar A=1$, $\bar B=1$, $\bar Q=-1$, $\bar R=0$, $\bar S=0$, such that the RDARE reduces to $\bar P+1=0$, which corresponds to $\bar K=1$, such that $\bar A-\bar B \bar K=0$ is stable. However, $R+B^\top \bar P B = 0$, and the stabilizing solution of the RDARE does not solve the DARE. $\hfill \square$
\end{Example}

We prove next that, assuming that $(A,B)$ is controllable, quadratic strict $(x,u)$-pre-dissipativity is not only necessary and sufficient for the existence of the stabilizing solution of the DARE (a proof can be found in~\cite[Corollary~8]{Zanon2016b}), but it also entails existence of the stabilizing solution to the RDARE.
\begin{Lemma}
	\label{lem:rdare_existence}
	Assume that 
	$(A,B)$ is controllable, and quadratic strict $(x,u)$-pre-dissipativity holds. Then the RDARE does have a stabilizing solution.
\end{Lemma}
\begin{IEEEproof}
	Without loss of generality we can assume that $A$ is invertible. Otherwise, since $(A,B)$ is controllable, we choose a feedback $\hat K$ such that $A-B\hat K$ is invertible. 
	Then the pre-stabilized system is still controllable 
	and by Lemma \ref{lem:pre_stabilized_system} it is still pre-dissipative and the corresponding DARE (and thus also the RDARE) has the same solution as for the original system.
	
	Assuming full rank of $A$, we first observe that the original DARE and the rotated DARE deliver solution matrices satisfying $P_\Lambda = P+\Lambda$. Consequently, we can assume without loss of generality that quadratic strict $(x,u)$-pre-dissipativity holds for $\Lambda=0$, i.e., $H\succ0$. Then we observe that
	\begin{align*}
		&\matr{ll}{-\bar A^\top & 0 \\ -\bar B^\top  & I }\matr{ll}{Q & S^\top \\ S & R}\matr{rr}{-\bar A & -\bar B \\ 0 & I } \\
		= \ & \matr{ll}{-\bar A^\top Q & -\bar A^\top S^\top \\ -\bar B^\top Q + S & -\bar B^\top S^\top + R}\matr{rr}{-\bar A & -\bar B \\ 0 & I } \\
		= \ & \matr{ll}{\bar A^\top Q \bar A & \bar A^\top Q \bar B - \bar A^\top S^\top \\ \bar B^\top Q \bar A - S \bar A & \bar B^\top Q \bar B - S \bar B - \bar B^\top S^\top + R} \\
		= \ & - \matr{ll}{\bar Q & \bar S^\top \\ \bar S & \bar R},
	\end{align*}
	such that, since the cost matrix $H$ is pre- and post-multiplied by a full-rank matrix, it holds that $\bar H \prec 0$. We observe that controllability of $(A,B)$ implies controllability of $(\bar A, \bar B)$. 
	Moreover, for any DARE / RDARE, changing the sign of the stage cost $H$ or $\bar H$ entails changing the sign of the solution $P$ or $\bar P$, while $K$ and $\bar K$ remain unchanged. 
	Consequently, the RDARE does have a stabilizing solution $\bar P_\mathrm{s}\prec0$. 
\end{IEEEproof}

As anticipated above, however, the fact that the RDARE does have a stabilizing solution does not entail that the DARE has an antistabilizing solution. We provide a sufficient condition in the next theorem.

\begin{theorem}
	\label{thm:antistabilizing_existence}
	Assume that $(A,B)$ is controllable and quadratic strict $(x,u)$-pre-dissipativity holds. Select $\hat K$ such that $A_{\hat K}:=A-B\hat K$ is invertible. 
	Define $D:=R-S_{\hat K}A_{\hat K}^{-1}B$, with $S_{\hat K}=S-R\hat K$. Assume that $D$ is invertible. Then, the antistabilizing solution of the DARE exists.
\end{theorem}

Before proving the theorem, let us briefly comment on the implications of pre-stabilizing the system.
\begin{Remark}
	By Lemma~\ref{lem:rotated_equivalence} and~\cite[Lemma~2]{Zanon2014d}, since quadratic strict $(x,u)$-pre-dissipativity holds, we can assume without loss of generality $H\succ0$ and, hence, $R\succ0$, i.e., it is invertible. Consequently, the nonsingularity of $R-SA^{-1}B$ is equivalent to the nonsingularity of
	\begin{align}
		\label{eq:nonsingular_matrix}
		\matr{ll}{ R & S \\ B & A}.
	\end{align}
	For the same reason, the nonsingularity of $R-S_{\hat K}A_{\hat K}^{-1}B$ is equivalent to the nonsingularity of
	\begin{align}
		\label{eq:nonsingular_matrixK}
		\matr{ll}{ R & S - R\hat K \\ B & A-B\hat K} = \matr{ll}{ R & S \\ B & A}\matr{ll}{I & -\hat K \\ 0  & I }.
	\end{align}
	Since the last matrix is full rank, nonsingularity of $R-SA^{-1}B$ is equivalent to nonsingularity of $R-S_{\hat K}A_{\hat K}^{-1}B$, i.e., that property is independent of any pre-stabilization of the system, but allows us to characterize the existence of the minimal solution also in case $A$ is not invertible.
	Finally, the condition $A-BR^{-1}S=A_{\hat K}-BR^{-1}S_{\hat K}$ nonsingular is another equivalent condition which is sometimes found in the literature, see, e.g.,~\cite{Ran1993a}.
\end{Remark}

\begin{IEEEproof}[Proof of Theorem~\ref{thm:antistabilizing_existence}]
	We first prove that the claim can be verified by checking it for a pre-stabilized system instead of the original one. This is particularly important in case $A$ is not invertible. Controllability ensures that there exist a feedback matrix $\hat K$ such that the eigenvalues of $A-B\hat K$ can be chosen arbitrarily. In particular, this entails that $A-B\hat K$ can be made invertible.
	
	By~\cite[Lemma~1]{Zanon2014d}, we know that the DARE formulated with system matrices $(A_{\hat K},B)$ and cost matrix $H_{\hat K}$ defined as per~\eqref{eq:pre_stabilized_cost} yields the same matrix as the original DARE, i.e., $P_{\hat K}=P$, while the feedback matrix is given by $K_{\hat K}=K-\hat K$. This entails that the solutions of the pre-stabilized DARE coincide with those of the original DARE.

	We can therefore exploit Lemma~\ref{lem:rdare_existence} to conclude existence of a stabilizing solution to the RDARE. Moreover, because $R-S_{\hat K}A_{\hat K}^{-1}B$ is nonsingular the stabilizing solution to the RDARE coincides with the antistabilizing solution to the DARE~\cite[Proposition~2]{Ionescu1996}.
\end{IEEEproof}

We establish next a relation between $\bar P_\mathrm{s}$ or $P_\mathrm{a}$ and cost rotations.

\subsection{Positive-Semidefinite Cost Rotations}

	Next, we need to establish one additional useful result. To that end, we first prove the following lemma.
	\begin{Lemma}
		Consider a DARE and rotate the stage cost matrix using matrix $\Lambda$. The corresponding RDARE coincides with the RDARE associated with the original cost rotated using matrix $\Lambda$.
	\end{Lemma}
	\begin{IEEEproof}
		We observe that the rotated cost yields the following stage cost matrices for the RDARE:
		\begin{align*}
			\bar Q_\lambda &= -\bar A^\top Q_\lambda \bar A, \\
			&= -\bar A^\top (Q +\Lambda - A^\top\Lambda A) \bar A \\
			&= \bar Q -\bar A^\top \Lambda \bar A + \Lambda, 
			\\
			\bar S_\lambda &= S_\lambda \bar A - \bar B^\top Q_\lambda \bar A, \\
			&= (S - B^\top \Lambda A) \bar A - \bar B^\top (Q +\Lambda - A^\top\Lambda A) \bar A, 
			\\
			&= S\bar A - B^\top \Lambda - \bar B^\top Q \bar A - \bar B^\top\Lambda \bar A + \underbrace{\bar B^\top A^\top}_{=B^\top}\Lambda, 
			\\
			&= \bar S - \bar B^\top \Lambda \bar A, \\
			\bar R_\lambda &= -R_\lambda + S_\lambda \bar B + \bar B^\top S_\lambda^\top - \bar B^\top Q_\lambda \bar B \\
			&= -(R - B^\top\Lambda B) + (S - B^\top \Lambda A) \bar B \\
			&\hspace{11pt} + \bar B^\top (S^\top - A^\top \Lambda B) - \bar B^\top (Q +\Lambda - A^\top\Lambda A)\bar B\\
			&= -R + B^\top\Lambda B + S\bar B - B^\top \Lambda \underbrace{A \bar B}_{=B} \\
			&\hspace{11pt} + \bar B^\top S^\top - \underbrace{\bar B^\top A^\top}_{=B^\top} \Lambda B - \bar B^\top Q \bar B - \bar B^\top \Lambda \bar B \\
			&\hspace{11pt} + \underbrace{\bar B^\top A^\top}_{=B^\top} \Lambda \underbrace{A \bar B}_{=B} \\
			&= \bar R - \bar B^\top \Lambda \bar B.
		\end{align*}
		Consequently, we obtain
		\begin{align*}
			\bar H_\Lambda = \matr{ll}{\bar Q + \Lambda - \bar A^\top \Lambda \bar A & \bar S^\top - \bar A^\top \Lambda \bar B \\ \bar S - \bar B^\top \Lambda \bar A & \bar R - \bar B^\top \Lambda \bar B},
		\end{align*}
		which, by definition, is the stage cost matrix obtained by rotating the original RDARE cost with matrix $\Lambda$.
	\end{IEEEproof}

	We are now ready to prove the following.
	\begin{theorem}
		\label{thm:rotation_barPs}
		Assume that $(A,B)$ is controllable and quadratic strict $(x,u)$-pre-dissipativity holds. Rotate the cost using $\Lambda=-\bar P_\mathrm{s}$. Then the rotated stage cost matrix $H_\Lambda$ is positive semidefinite, i.e., $H_\Lambda\succeq0$.
	\end{theorem}
	\begin{IEEEproof}
		In order to obtain the proof we first focus on the stage cost matrix of the RDARE, which reads
		\begin{align*}
			\bar H_\Lambda &:= \matr{ll}{\bar Q - \bar P_\mathrm{s} + \bar A^\top \bar P_\mathrm{s} \bar A & \bar S^\top + \bar A^\top \bar P_\mathrm{s} \bar B \\ \bar S + \bar B^\top \bar P_\mathrm{s} \bar A & \bar R + \bar B^\top \bar P_\mathrm{s} \bar B}.
		\end{align*}
		By Lemma~\ref{lem:rdare_existence} we know that if $H\succ0$ the stage cost of the RDARE is negative definite, such that $\bar R\prec 0$, $\bar P_\mathrm{s} \prec 0$. This entails
		\begin{align*}
			\bar R_\lambda = \bar R + \bar B^\top \bar P_\mathrm{s} B \prec 0.
		\end{align*}
		Then we observe that the Schur complement of $\bar R$ in $\bar H_\Lambda$ yields the RDARE, which is solved by $\bar P_\mathrm{s}$. Consequently, $\bar H_\Lambda \preceq 0$. 
		
		Finally, we recall that 
		\begin{align*}
			-\bar H_\Lambda = &\matr{ll}{-\bar A^\top & 0 \\ -\bar B^\top  & I } H_\Lambda \matr{rr}{-\bar A & -\bar B \\ 0 & I } \succeq 0.
		\end{align*}
		Since in the rhs of the equation above $H_\Lambda$ is pre- and post-multiplied by a full rank matrix and its transpose, we conclude that $H_\Lambda\succeq0$.
	\end{IEEEproof}
	Note that, as a result of Lemma~\ref{lem:rotated_equivalence}, after rotating with $\Lambda=-\bar P_\mathrm{s}$ we have that the new solution is $\bar P_\mathrm{s}=0$, consequently, if it exists, $P_\mathrm{a}=0$.
		
	Finally, we discuss next further useful reformulations and properties.

\subsection{Properties of the Cost-to-Go and DARE Solutions}

Finite-horizon LQR problems are characterized by the Riccati iterations, defined by~\eqref{eq:RE}, which
is fully equivalent to
\begin{align*}
	P_{n+1} =\ &  A_{K_{n+1}}^\top P_n A_{K_{n+1}}  \\
	&+  \underbrace{Q - S^\top K_{n+1} - K_{n+1}^\top S + K_{n+1}^\top R K_{n+1}}_{=:Q_{K_{n+1}}},
\end{align*}
where $K_{n+1}$ is defined as per~\eqref{eq:RE:K}.

%
%
%

\begin{Lemma}
	\label{lem:pd_stability}
	Assume that $P_{n+1}\succeq P_n$. Then, $P_{n+2}\succeq P_{n+1}$. 
\end{Lemma}
\begin{IEEEproof}
	We observe that
	\begin{align*}
		P_{n+2} &= Q_{K_{n+2}} + A_{K_{n+2}}^\top P_{n+1} A_{K_{n+2}} \\
		&\succeq Q_{K_{n+2}} + A_{K_{n+2}}^\top P_{n} A_{K_{n+2}} \\
		&\succeq Q_{K_{n+1}} + A_{K_{n+1}}^\top P_{n} A_{K_{n+1}} = P_{n+1},
	\end{align*}
	where we used optimality of the one-step-ahead problem with terminal cost matrix $P_n$ to obtain the second inequality, while the first inequality stems from the assumption $P_{n+1}\succeq P_n$. 
	%
\end{IEEEproof}
This lemma entails that, if the LQR problem over a horizon $n=1$ yields a cost-to-go which is no smaller than the terminal cost, then the sequence of matrices $P_n$ yielded by~\eqref{eq:RE} is monotonic.

Let us denote the set of solutions to the DARE as $\mathcal{P}$ and the set of the non stabilizing solutions as $\mathcal{\bar P}:=\mathcal{P}\setminus \{P_\mathrm{s}\}$. 
We call a solution $P_\mathrm{a}\in \mathcal{\bar P}$ {\em antistabilizing}, if for the corresponding feedback matrix as $K_\mathrm{a}$ all eigenvalues $\mu$ of $A_{\mathrm{a}}=A-BK_\mathrm{a}$ satisfy $|\mu|>1$.
\begin{Lemma}[{\cite[Lemmas 3.1--3.2]{Ran1993a}}]
	\label{lem:delta_riccati}
	Suppose that $P$ is an arbitrary solution to the DARE. Provided the stabilizing and the antistabilizing solution to this DARE exists, define 
	\begin{align}
		R_\mathrm{s} := R + B^\top P_\mathrm{s} B && R_\mathrm{a} := R + B^\top P_\mathrm{a} B.
	\end{align}
	Then, $\Delta_\mathrm{a}:=P-P_\mathrm{a}$ and $\Delta_\mathrm{s}:=P-P_\mathrm{s}$ satisfy the algebraic Riccati equations
	\begin{align}
		\Delta_\mathrm{a} &= A_\mathrm{a}^\top \Delta_\mathrm{a} A_\mathrm{a} -A_\mathrm{a}^\top \Delta_\mathrm{a} B (R_\mathrm{a} + B^\top \Delta_\mathrm{a} B)^{-1} B^\top \Delta_\mathrm{a} A_\mathrm{a}, \label{eq:dare_diff_d} \\
		\Delta_\mathrm{s} &= A_\mathrm{s}^\top \Delta_\mathrm{s} A_\mathrm{s} -A_\mathrm{s}^\top \Delta_\mathrm{s} B (R_\mathrm{s} + B^\top \Delta_\mathrm{s} B)^{-1} B^\top \Delta_\mathrm{s} A_\mathrm{s}. \label{eq:dare_diff_s}
	\end{align}
\end{Lemma}

\begin{Lemma}
	\label{lem:delta_riccati_iterations}
	Suppose that $P_+,P$ solve the Riccati iteration 
	\begin{subequations}
		\begin{align}
			P_{+} =\ & Q + A^\top P A - (S^\top + A^\top P B)K, \\
			K = \ & (R+B^\top P B)^{-1} (S+B^\top P A).
		\end{align}
	\end{subequations}
	Provided the stabilizing and the antistabilizing solution to this DARE exist, define 
	\begin{align}
		R_\mathrm{s} := R + B^\top P_\mathrm{s} B && R_\mathrm{a} := R + B^\top P_\mathrm{a} B.
	\end{align}
	Then, $\Delta_\mathrm{a}:=P-P_\mathrm{s}$, $\Delta_\mathrm{a}^+:=P_+-P_\mathrm{s}$ and $\Delta_\mathrm{a}:=P-P_\mathrm{a}$, $\Delta_\mathrm{a}^+:=P_+-P_\mathrm{a}$ satisfy the algebraic Riccati equations
	\begin{align}
		\Delta_\mathrm{a}^+ &= A_\mathrm{a}^\top \Delta_\mathrm{a} A_\mathrm{a} -A_\mathrm{a}^\top \Delta_\mathrm{a} B (R_\mathrm{a} + B^\top \Delta_\mathrm{a} B)^{-1} B^\top \Delta_\mathrm{a} A_\mathrm{a}, \label{eq:daiter}\\
		\Delta_\mathrm{s}^+ &= A_\mathrm{s}^\top \Delta_\mathrm{s} A_\mathrm{s} -A_\mathrm{s}^\top \Delta_\mathrm{s} B (R_\mathrm{s} + B^\top \Delta_\mathrm{s} B)^{-1} B^\top \Delta_\mathrm{s} A_\mathrm{s}.
	\end{align}
\end{Lemma}
\begin{IEEEproof}
	The proof follows along the same lines as Lemma~\ref{lem:delta_riccati} and~\cite[Lemmas 3.1--3.2]{Ran1993a}.
\end{IEEEproof}

\begin{Remark}
	Note that these lemmas entail that the differences between solutions of a DARE can be interpreted as solutions to an optimal control problem which only penalizes the input with matrix $R_\mathrm{a}$ or $R_\mathrm{s}$, respectively. Consequently, we can apply Theorem~\ref{thm:stability} to prove some useful properties of all solutions to the DARE.
\end{Remark}

We recall next that $P_\mathrm{s}\succeq P \succeq P_\mathrm{a}$. A proof is provided in~\cite{Ran1993a}, but we provide an alternative and shorter one next.
\begin{Lemma}
	\label{lem:riccati_solutions_order_s}
	Whenever $P_\mathrm{s}$ exists, it holds that $P_\mathrm{s}\succeq P$, for all symmetric $P$ solving the DARE. 
\end{Lemma}
\begin{IEEEproof}
	%
	The claim can equivalently be formulated as $\Delta_\mathrm{s}\preceq0$. We observe that, by Theorem~\ref{thm:dare_pd} we have
	\begin{align*}
		R_\mathrm{s} + B^\top \Delta_\mathrm{s} B = R + B^\top P B \succ 0.
	\end{align*}
	This entails that
	\begin{align*}
		M_\mathrm{s} := A_\mathrm{s}^\top \Delta_\mathrm{s} B (R_\mathrm{s} + B^\top \Delta_\mathrm{s} B)^{-1} B^\top \Delta_\mathrm{s} A_\mathrm{s} \succeq 0.
	\end{align*}
	We observe that~\eqref{eq:dare_diff_s} also reads
	\begin{align*}
		\Delta_\mathrm{s} &= A_\mathrm{s}^\top \Delta_\mathrm{s} A_\mathrm{s} - M_\mathrm{s},
	\end{align*}
	such that, since $M_\mathrm{s} \succeq 0$, this directly entails 
	\begin{align*}
		\Delta_\mathrm{s} &= -\sum_{k=0}^\infty \underbrace{A_\mathrm{s}^\top M_\mathrm{s} A_\mathrm{s}}_{\succeq0}\preceq0.
	\end{align*}
\end{IEEEproof}

\begin{Lemma}
	\label{lem:riccati_solutions_order_d}
	Whenever $P_\mathrm{a}$ exists, it holds that $P_\mathrm{a}\preceq P$, for all symmetric $P$ solving the DARE. 
\end{Lemma}
\begin{IEEEproof}
	The claim can equivalently be formulated as $\Delta_\mathrm{a}\succeq0$. In order to prove the result, we first observe that $A_\mathrm{a}$ has all its eigenvalues outside the unit circle, such that $A_\mathrm{a}^{-1}$ exists and has all its eigenvalues inside the unit circle. Then,~\eqref{eq:dare_diff_d} can be rewritten as
	\begin{align*}
		A_\mathrm{a}^{-\top} \Delta_\mathrm{a} A_\mathrm{a}^{-1} &= \Delta_\mathrm{a}  - \Delta_\mathrm{a} B (R_\mathrm{a} + B^\top \Delta_\mathrm{a} B)^{-1} B^\top \Delta_\mathrm{a}.
	\end{align*}
	By rearranging terms we have
	\begin{align*}
		\Delta_\mathrm{a} = A_\mathrm{a}^{-\top} \Delta_\mathrm{a} A_\mathrm{a}^{-1} +  \Delta_\mathrm{a} B (R_\mathrm{a} + B^\top \Delta_\mathrm{a} B)^{-1} B^\top \Delta_\mathrm{a}.
	\end{align*}
	
	By Theorem~\ref{thm:dare_pd} we have
	\begin{align*}
		R_\mathrm{a} + B^\top \Delta_\mathrm{a} B = R + B^\top P B \succ 0.
	\end{align*}
	This entails that
	\begin{align*}
		M_\mathrm{a} := \Delta_\mathrm{a} B (R_\mathrm{a} + B^\top \Delta_\mathrm{a} B)^{-1} B^\top \Delta_\mathrm{a} \succeq 0.
	\end{align*}
	The remainder of the proof follows along the same lines as for Lemma~\ref{lem:riccati_solutions_order_s}.
\end{IEEEproof}


\begin{Lemma}
	Assume that both $P_\mathrm{s}$ and $P_\mathrm{a}$ exist  and define $\Xi_\mathrm{s}:=P_\mathrm{s}-P_\mathrm{a}$. Then $\Xi_\mathrm{s} \succ 0$, i.e., $P_\mathrm{s}\succ P_\mathrm{a}$.
\end{Lemma}
\begin{IEEEproof}
	Lemma~\ref{lem:riccati_solutions_order_d} already established $P_\mathrm{s}\succeq P_\mathrm{a}$, so we are left with proving that the inequality is strict. 
	
	In order to prove this result, it is convenient to use the reformulation proposed in Lemmas~\ref{lem:delta_riccati}--\ref{lem:delta_riccati_iterations}, i.e., use $\Delta_\mathrm{a}$ instead of $P$, and observe that $\Xi_\mathrm{s}$ can be interpreted as the solution of an infinite-horizon problem with a positive semidefinite stage cost and a suitably defined terminal cost/constraint for a linear system defined by matrices $(A_\mathrm{a},B)$. In particular, the terminal cost related to $P_\mathrm{a}$ is $0$. Since $A_\mathrm{a}$ is unstable and the feedback corresponding to $\Xi_\mathrm{s}$ is stabilizing, this entails that, for any initial state $x\neq0$, the optimal input cannot be $0$ at all times, which, together with $R_\mathrm{a}\succ0$, entails that the optimal cost-to-go needs to be strictly positive. In turn, this entails $\Xi_\mathrm{s}\succ0.$
\end{IEEEproof}

\section{Stabilizing Terminal Costs}
\label{sec:main}

We are now ready to prove the main results. We proceed in three steps in the next three sections: We first establish asymptotic stability of the optimal solutions of the infinite-horizon LQR problem assuming existence of the antistabilizing solution $P_\mathrm{a}$, then extend this result to the case in which $P_\mathrm{a}$ does not exist and finally address the receding horizon feedback law.

\subsection{The Antistabilizing Solution Does Exist}

Before discussing the case of indefinite stage costs, let us first consider the simpler case of positive semidefinite costs. To that end, we first need to prove the following lemma, which establishes some form of monotonicity of the sequence of cost-to-go matrices $P_n$.
\begin{Lemma}
	\label{lem:sequence_nondecreasing_simple}
	Consider an LQR problem formulated with terminal cost matrix $P_0=P_\mathrm{a} + \alpha \Xi_\mathrm{s}$, with scalar $0< \alpha\ < 1$ and $\Xi_\mathrm{s}:=P_\mathrm{s}-P_\mathrm{a}$. Then, $P_1\neq P_0$ and $P_1\succeq P_0$.
\end{Lemma}
\begin{IEEEproof}
	In order to prove the result it will be convenient to look at the solutions $\Delta_n$ of the Riccati iteration \eqref{eq:daiter} for $\Delta_{\mathrm{a}}$, which is initialized with terminal cost matrix $\Delta_0=\alpha \Xi_\mathrm{s}$. The claim then reads $\Delta_1\neq\Delta_0$, $\Delta_1\succeq\Delta_0$.
	
	We first prove the second claim. We use Lemma~\ref{lem:delta_riccati_iterations} and the Woodbury matrix identity to obtain
	\begin{align*}
		\Delta_1 &= A_\mathrm{a}^\top \Delta_0 A_\mathrm{a} -A_\mathrm{a}^\top \Delta_0 B (R_\mathrm{a} + B^\top \Delta_0 B)^{-1} B^\top \Delta_0 A_\mathrm{a}, \\
		&= A_\mathrm{a}^\top (\Delta_0^{-1}+B R_\mathrm{a}^{-1} B^\top)^{-1} A_\mathrm{a}  \\
		&= A_\mathrm{a}^\top (\alpha^{-1} \Xi_\mathrm{s}^{-1}+B R_\mathrm{a}^{-1} B^\top)^{-1} A_\mathrm{a} \\
		&= \alpha A_\mathrm{a}^\top (\Xi_\mathrm{s}^{-1}+ \alpha B R_\mathrm{a}^{-1} B^\top)^{-1} A_\mathrm{a}.
	\end{align*}
	We observe that, by Lemma~\ref{lem:delta_riccati} and the Woodbury matrix identity we have
	\begin{align*}
		\alpha \Xi_\mathrm{s} &= \alpha \left (A_\mathrm{a}^\top \Xi_\mathrm{s} A_\mathrm{a} - A_\mathrm{a}^\top \Xi_\mathrm{s} B ( R_\mathrm{a} + B^\top \Xi_\mathrm{s} B)^{-1} B^\top \Xi_\mathrm{s} A_\mathrm{a} \right ), \\
		&= \alpha A_\mathrm{a}^\top (\Xi_\mathrm{s}^{-1}+B R_\mathrm{a}^{-1} B^\top)^{-1} A_\mathrm{a}.
	\end{align*}
	Because $\Xi_\mathrm{s}^{-1}+B R_\mathrm{a}^{-1} B^\top \succeq \Xi_\mathrm{s}^{-1}+\alpha B R_\mathrm{a}^{-1} B^\top \succ 0$, then we have $(\Xi_\mathrm{s}^{-1}+\alpha B R_\mathrm{a}^{-1} B^\top)^{-1} \succeq (\Xi_\mathrm{s}^{-1}+B R_\mathrm{a}^{-1} B^\top)^{-1} \succ 0$, which entails $\Delta_1 \succeq \Delta_0$, i.e., $P_1\succeq P_0$.
	
	The first claim is directly obtained by contradiction, as $\Delta_1=\Delta_0$ implies that $B R_\mathrm{a}^{-1} B^\top = \Xi_\mathrm{s}^{-1}+\alpha B R_\mathrm{a}^{-1} B^\top$, which is clearly impossible, since by Theorem~\ref{thm:dare_pd} we have $R_\mathrm{a}\succ0$.
%
\end{IEEEproof}

\begin{theorem}
	\label{thm:stability}
	Consider Problem~\eqref{eq:lqr_terminal_cost}
	with $H\succeq0$, 
	and $P^\mathrm{f}\succ0$ and let $(A,B)$ be stabilizable. 
	Then, if quadratic strict $(x,u)$-pre-dissipativity holds, the DARE has a stabilizing solution and the corresponding optimal value function and feedback coincide with those of the infinite horizon problem~\eqref{eq:lqr_terminal_cost}. $\hfill\square$
\end{theorem}
\begin{IEEEproof}
		Since quadratic strict $(x,u)$-pre-dissipativity holds the rotated cost is positive definite. Together with the stabilizability of $(A,B)$ this implies that the rotated DARE has a stabilizing solution \cite[Section 3.3]{Anderson1990} 
		which by Lemma~\ref{lem:pre_stabilized_system} has a corresponding unique stabilizing solution of the DARE for the original system. This shows the first claim.
	
	In order to prove the second claim, proceeding backwards in time starting from the final time, we write the dynamic programming recursion
	\begin{align*}
		W_{n+1}(x) = \min_v \ \ell(x,v) + W_{n}( A x + B v),
	\end{align*}
	with $W_0(x)=x^\top P^\mathrm{f} x$
. We observe that, due to the quadratic form of $\ell$, for all time instants $n$ 
the functions $W_n(x)$ are quadratic, i.e., of the form $W_n(x)=x^\top P_n x$ for symmetric and positive definite matrices $P_n$. Since the pair $(A,B)$ is stabilizable, the $W_n$ are also uniformly bounded on compact sets. We prove next that the sequence $P_n$ converges. By construction $P_n\succeq0$. We consider first the case $P^\mathrm{f}\prec P_\mathrm{s}$, which entails $P_n \prec P_\mathrm{s}$ for all $n<\infty$. Moreover, as discussed in Lemma~\ref{lem:rdare_existence}, we have $P_\mathrm{a}=\bar P_\mathrm{s} \preceq 0$, such that $P^\mathrm{f} \succ P_\mathrm{a}$. Lemma~\ref{lem:sequence_nondecreasing_simple} then entails that $P_n$ is bounded from below by the sequence $\hat P_n$ obtained using terminal cost $P_\mathrm{a}+\alpha\Xi_\mathrm{s} \preceq P^\mathrm{f}$. Note that such terminal cost always exists for a sufficiently small $\alpha>0$, since $P^\mathrm{f} \succ P_\mathrm{a}$. The sequence $\hat P_n$ is nondecreasing and bounded, such that it must converge to a fixed value. The case $P^\mathrm{f}\nprec P_\mathrm{s}$ is obtained as a direct extension to the argument above. Indeed, we know that for any $\tilde P_0\succ P_\mathrm{s}$ we obtain a sequence $\tilde P_{n+1}\preceq \tilde P_n$ which converges to $P_\mathrm{s}$, see, e.g.,~\cite{Gruene2017,Rawlings2017}. This immediately provides an upper bound for the case $P^\mathrm{f}\nprec P_\mathrm{s}$, while the lower bound remains unaltered and we obtain once more that the sequence $P_n$ converges to a fixed value. Because we will prove next that $\lim_{n\to\infty}\hat P_n=P_\mathrm{s}$, this will also entail that $\lim_{n\to\infty}P_n=P_\mathrm{s}$. Consequently, we will assume without loss of generality that $P^\mathrm{f}=P_\mathrm{a}+\alpha\Xi_\mathrm{s}$, such that $P_n=\hat P_n$.

Since the sequence converges, using dynamic programming we obtain that there exists a matrix $K$ such that the infinite horizon optimal control is given by the feedback law $u_k^\star = -Kx_k^\star$, i.e., by $F(x)=-Kx$. In order to prove the claim, we need to show that all eigenvalues $A-BK$ are inside the unit circle, i.e., that $P_\infty$ is the unique stabilizing solution of the DARE.
	
	We observe that either the optimal linear feedback policy $F(x)=-Kx$ yields a closed-loop matrix $A-BK$ with stable eigenvalues (inside or on the boundary of the unit circle with those on the boundary being semi-simple); or one or more eigenvalues are outside the unit circle or at the boundary and not semi-simple. 
	Because any feedback policy which does not stabilize the system yields a state trajectory for which at least one component of the state diverges to $\pm\infty$, all policies of the latter type incur an unbounded terminal cost as $N\to\infty$, which contradicts the fact that $V(x) = x^\top P_\infty x$ is finite for all $x$. 
	
	

Since this excludes the case of diverging solutions, it remains to exclude the possibility that the optimal closed-loop matrix has one or more eigenvalues on the unit circle, i.e., $|\mu|=1$ for some eigenvalue $\mu$ of $A_K:=A-BK$.

	To that end, 
	consider the reformulation of problem~\eqref{eq:lqr_terminal_cost} given in~\eqref{eq:lqr_terminal_cost2_rot}, where we apply the change of variable $\bar u_k = u_k-Kx_k$.
	Since $-K$ is the infinite-horizon optimal feedback law, the optimal control of the rotated and pre-stabilized problem~\eqref{eq:lqr_terminal_cost2_rot} is $\bar u_k=0$. 
	Consequently, the optimal cost is given by
	\begin{align}
		V(x_0) = \lim_{N\to\infty} &\sum_{k=0}^{N-1} (A_K^kx_0)^\top (Q_K+\Lambda-A_K^\top \Lambda A_K) A_K^kx_0 \nonumber \\
		&+ x_N^\top (P^\mathrm{f}+\Lambda) x_N.\label{eq:Vsum}
	\end{align}
	As we proved before, $\lim_{N\to \infty}\|x_N\|< \infty$, such that the terminal cost is bounded for $N\to\infty$.
	
	In order to prove that $|\mu|<1$ for all eigenvalues of $A_K$, let us proceed by contradiction: let us assume that there is an eigenvalue with $|\mu|=1$, implying $\|A_K^kv_\mu\| = \|v_\mu\|>0$ for the corresponding eigenvector $v_\mu$. Since quadratic strict $(x,u)$-pre-dissipativity holds, the matrix $Q_K+\Lambda-A_K^\top \Lambda A_K$ is positive definite, implying the existence of a constant $c>0$ such that 
	\[ (v_\mu A_K^k)^\top (Q_K+\Lambda-A_K^\top \Lambda A_K) A_K^kv_\mu \ge c\|A_K^kv_\mu \|^2 \ge c \|v_\mu\|^2.\]
This implies that the sum in \eqref{eq:Vsum} diverges for $x_0=v_\mu$, which again contradicts the finiteness of the optimal cost.
%
\end{IEEEproof}

We are now ready to extend the result above to indefinite stage costs.

\begin{theorem}
	\label{thm:stability_general}
	Consider a strictly $(x,u)$-pre-dissipative LQR problem, and assume that the antistabilizing  solution $P_\mathrm{a}$ of the corresponding DARE exists. Then, for any terminal cost matrix $P_0=P_\mathrm{a} + E$, with $E \succ 0$, the optimal solution is given by the unique stabilizing solution of the corresponding DARE. In particular, the origin is an asymptotically stable equilibrium of the optimally controlled system.
\end{theorem}
\begin{IEEEproof}
		In order to prove the result it will be convenient to look at the matrix $\Delta_\mathrm{a}$, which defines the optimal cost for the infinite-horizon optimal control problem with zero penalization of the state in the stage cost 
		and terminal cost matrix $\Delta_0=E \succ0$.
		
		By Theorem~\ref{thm:stability} we have that, the terminal cost being positive definite and the stage cost positive semidefinite, the solution to the infinite-horizon LQR must coincide with the unique stabilizing solution of the DARE. 
\end{IEEEproof}

\begin{Corollary}
	\label{cor:P_order}
	Consider any $P\in\mathcal{\bar P}$, i.e., $P\neq P_\mathrm{s}$ solving the DARE~\eqref{eq:DARE}. Then, 
	if $P_\mathrm{a}$ exists, $P \nsucc P_\mathrm{a}$.
\end{Corollary}
\begin{IEEEproof}
	Our proof will exploit uniqueness of the stabilizing solution, i.e., that $P\neq P_\mathrm{s}$ cannot be a stabilizing solution. 
	
	
	Assume that $P\in\mathcal{\bar P}$ and $P_\mathrm{s}\neq P \succ P_\mathrm{a}$. Because $P$ is a symmetric solution to the DARE and $R+B^\top P B \succ 0$, when using $P$ as terminal cost matrix, the cost-to-go matrix of the LQR associated with the DARE is $P$ for any horizon length. 
	However, by Theorem~\ref{thm:stability_general}, because $P \succ P_\mathrm{a}$, the infinite-horizon LQR yields the unique stabilizing solution to the DARE~\eqref{eq:DARE}. Since this contradicts uniqueness of the stabilizing solution, $P \nsucc P_\mathrm{a}$.
\end{IEEEproof}
\begin{Remark}
	Note that the results above implicitly prove that the unique stabilizing solution is the only solution to the DARE satisfying $P\succ P_\mathrm{a}$ and all other solutions can only satisfy $P\succeq P_\mathrm{a}$.
\end{Remark}
\begin{Remark}
	Note that, though we proved the result for terminal cost matrix $P_0=P_\mathrm{a} + E$, the result clearly holds for terminal cost matrix $P_0=\bar P + E$, with $\bar P$ any solution of the DARE. However, referring the result to $P_\mathrm{a}$ yields the least restrictive condition as, by Lemma~\ref{lem:riccati_solutions_order_d} we have $P_\mathrm{a}\preceq\bar P$. Moreover, this sufficient condition is also close to being necessary necessary, as selecting $E=0$ yields a cost-to-go matrix $P_n=P_\mathrm{a}$ for all $n$, with a corresponding destabilizing feedback. A similar reasoning holds for $E=P-P_\mathrm{a}$, where $P \neq P_\mathrm{s}$ is any other non-stabilizing symmetric solution to the DARE, in which case the closed-loop system has both stable and unstable eigenvalues. 
\end{Remark}

\subsection{The Antistabilizing Solution Does not Exist}

In this subsection we extend the previous results to also cover the case in which the antistabilizing solution of the DARE does not exist.

		\begin{theorem}
			\label{thm:convergence_no_antistab}
			Consider a strictly pre-dissipative LQR problem with controllable $(A,B)$ and let $\bar P_\mathrm{s}$ be the stabilizing solution of the RDARE (which exists according to Lemma \ref{lem:rdare_existence}). Choose the terminal cost matrix $P^\mathrm{f}=\bar P_\mathrm{s} + E$, with $E \succ 0$. Then the optimal solution is given by the unique stabilizing solution of the corresponding DARE. In particular, the origin is an asymptotically stable equilibrium of the optimally controlled system.

		\end{theorem}
		\begin{IEEEproof}
			We define $H^\eps$ as
			\begin{align*}
				0 \prec (1-\epsilon)H \preceq H^\eps := \matr{ll}{(1-\epsilon)Q & (1-\epsilon)S^\top \\ (1-\epsilon)S & (1-\epsilon)R+T^\eps} \preceq H,
			\end{align*}
			with $1>\epsilon>0$ sufficiently small and $T^\eps$ such that 
			$\epsilon R \succeq T^\eps \succ 0$. For this cost we have that the condition for the existence of the antistabilizing solution is the nonsingularity of
			\begin{align*}
				(1-\epsilon)(R - S_{\hat K} A_{\hat K}^{-1}B) + T^\eps.
			\end{align*}
			Because $T^\eps$ can be chosen in the cone $\epsilon R \succeq T^\eps \succ 0$, for all $\epsilon$ there exists $T^\eps$ such that the matrix above is nonsingular. To prove that, we write
				\begin{align*}
					NWN^{-1}=(1-\epsilon)(R - S_{\hat K} A_{\hat K}^{-1}B),
				\end{align*}
				in Jordan form and define 
				\begin{align*}
					T^\eps := N \bar W N^{-1},
				\end{align*}
				with $\bar W$ a diagonal matrix such that
				\begin{align*}
					\bar W_{ii} = \left \{ \begin{array}{ll}
						a & \text{if } W_{ii} = 0, \\
						0 & \text{otherwise}
					\end{array}\right.,
				\end{align*}
				with $0<a<\epsilon \sigma_{\mathrm{min}}(R)$. By construction, we then have that $(1-\epsilon)(R - S_{\hat K} A_{\hat K}^{-1}B) + T^\eps=  N(W+ \bar W) N^{-1}$ is nonsingular.
			Consequently, for all $\eps$ under consideration the antistabilizing solution $P_\mathrm{a}^\eps$ exists and matches the stabilizing solution $\bar P_\mathrm{s}^\eps$ of the RDARE. 
			Note that because $H^\eps \preceq H$, we have $(1-\epsilon)\bar P_\mathrm{s} \succeq P_\mathrm{a}^\eps =  \bar P_\mathrm{s}^\eps \succeq \bar P_\mathrm{s}$. Here $\bar P_\mathrm{s}$ denotes the stabilizing solution of the RDARE for $H$, for which $ P_\mathrm{s}^\eps \to \bar P_\mathrm{s}$ as $\epsilon \to 0$ holds. We recall that these are the RDARE stabilizing solutions, which decrease as the stage cost increases.
						
			One can now select $\epsilon$ small enough such that $E \succ P^\eps_\mathrm{a} - \bar P_\mathrm{s}$. Then, for all such $\eps$ by Theorem~\ref{thm:stability_general} the optimal solution for the problem with cost $H^\eps$ is given by the stabilizing solution $P_{\mathrm s}^\eps$ of the DARE for $H^\eps$, i.e., the optimal value function is given by 
			\[ V^\eps(x) = x^TP_{\mathrm s}^\eps x.\]
			Together with the ordering of the matrices $H^\eps$ this implies that 
			\[ (1-\eps) x^TP_{\mathrm s} x \le V^\eps(x) \le x^TP_{\mathrm s} x.\]
			Hence, for $\eps\to 0$ we get that $V(x) = x^TP_{\mathrm s} x$ which shows the claim. 
%
%
%
		\end{IEEEproof}
		
\begin{Remark}\label{rem:PLambda}
	Since by Lemma \ref{lem:rotated_equivalence} the optimal solutions for the original and the rotated cost coincide, under the condition of Theorem \ref{thm:convergence_no_antistab}, it follows that the origin is also exponentially stable for the infinite horizon optimal solutions for the rotated cost. This implies that according to \eqref{eq:Vsum} the optimal value function satisfies
	\begin{align*}
		V(x_0) & = x_0^TP_\lambda x_0 \\ & = \lim_{N\to\infty} \sum_{k=0}^{N-1} (A_K^kx_0)^\top (Q_K+\Lambda-A_K^\top \Lambda A_K) A_K^kx_0,
	\end{align*}
	as the terminal cost in \eqref{eq:Vsum} vanishes for $N\to\infty$. Hence, the matrix $P_\lambda$ is positive definite, since the first term in the above sum is positive definite and all others are positive semidefinite, because of quadratic strict $(x,u)$-pre-dissipativity. As Lemma \ref{lem:rotated_equivalence} states that 
	$P_{\lambda} = P_\mathrm{s} + \Lambda$, this implies that $P_\mathrm{s} + \Lambda \succ 0$.
\end{Remark}
	
	\subsection{Convergence Implies Exponential Stability}
	With the results above, we have proven that the cost-to-go matrix converges to the stabilizing solution to the DARE, such that exponential stability is obtained for an infinite horizon LQR with a properly selected terminal cost. We prove next that exponential stability is also obtained by the receding horizon feedback law for a sufficiently long but finite horizon $N$. 
	
		\begin{theorem}
			\label{thm:exponential_stability}
			Consider a strictly pre-dissipative LQR problem with controllable $(A,B)$ and let $\bar P_\mathrm{s}$ be the stabilizing solution of the RDARE (which exists according to Lemma \ref{lem:rdare_existence}). Select the terminal cost matrix as $P^\mathrm{f}\succ \bar P_\mathrm{s}$.
			Then, for any sufficiently large finite horizon $N$ the RH-OCP \eqref{eq:empc} yields a closed-loop system \eqref{eq:clsystem} for which the origin is globally asymptotically stable.
		\end{theorem}
		\begin{IEEEproof}
			By standard Lyapunov function arguments, exponential stability follows if we show that for all sufficiently large $N$ the iterate $P_{\lambda,N}$ of the rotated Riccati iteration is positive definite and satisfies
			\begin{align}
				\label{eq:stability_P}
				A_{K_N}^\top P_{\lambda,N} A_{K_N} - P_{\lambda,N} \prec 0. 
			\end{align}
			We know that 
			\begin{align*}
				P_{\lambda,n} 
				&= Q_{\lambda,K_{n}} + A_{K_{n}}^\top P_{\lambda,n-1} A_{K_{n}},
			\end{align*}
			with $Q_{\lambda,K_n} = Q_{K_n}+\Lambda-A_{K_n}^\top \Lambda A_{K_n}$. This implies
			\begin{align*}
				 &A_{K_{n}}^\top P_{\lambda,n} A_{K_{n}} - P_{\lambda,n}\\ &= - Q_{\lambda,K_{n}} + A_{K_{n}}^\top (P_{\lambda,n}-P_{\lambda,n-1}) A_{K_{n}},
			\end{align*}
			i.e., the decrease condition for stability \eqref{eq:stability_P} becomes
			\begin{align}
				\label{eq:stability_Q}
				Q_{\lambda,K_{n}} \succ A_{K_{n}}^\top (P_{\lambda,n}- P_{\lambda,n-1}) A_{K_{n}}.
			\end{align}
			Now, as $n\to\infty$, the left hand side converges to $Q_{\lambda,K}$, which is positive definite because ofquadratic strict $(x,u)$-pre-dissipativity, while the right hand side converges to $0$, since $K_n\to K$ and $P_{\lambda,n}\to P_{\mathrm{s}}+\Lambda$, which is positive definite as explained in Remark \ref{rem:PLambda}. Hence, there is $N>0$ such that \eqref{eq:stability_Q} holds and $P_{\lambda,n}\succ0$ for all $n\ge N$, which shows the claim.
%
%
	\end{IEEEproof}

\section{Relation With Known Results}
\label{sec:nonlinear}

In this section we connect our results to similar ones available in the literature.

\subsection{The Nonlinear Case}
In the companion paper~\cite{Gruene2024} we discuss the problem for the general nonlinear case. For the linear quadratic case the conditions derived therein require the existence of matrix $\Lambda$ such that
\begin{align}
	\label{eq:nl_conditions_lq}
	P^\mathrm{f} \succ -\Lambda, && H_\Lambda \succeq 0.
\end{align}
We observe that, by Theorem~\ref{thm:rotation_barPs} we have that the choice $\Lambda=-\bar P_\mathrm{s}$ yields a positive semidefinite rotated stage cost. If $P_\mathrm{a}$ does not exist, then using $-\Lambda \prec \bar P_\mathrm{s}$ entails that matrix $R-B^\top \Lambda B$ must have some negative eigenvalue, as $-B^\top \Lambda B \preceq B^\top \bar P_\mathrm{s} B$, and $-B^\top \Lambda B \neq B^\top \bar P_\mathrm{s} B$. 
In case $P_\mathrm{a}$ exists, using the same arguments as in Lemma~\ref{lem:sequence_nondecreasing_simple} one can prove that, for $P^\mathrm{f}=\bar P_\mathrm{s} - \alpha \Xi_\mathrm{s},$ with $\alpha > 0$ we have $P_1 \preceq P^\mathrm{f}$. The Schur complement of $H_\Lambda$ reads
\begin{align*}
	 &Q + \Lambda - A^\top \Lambda A \\
	 &\hspace{4em}- (S^\top - A^\top \Lambda B)(R - B^\top \Lambda B)^{-1} (S - B^\top \Lambda A) \\
	 &= P_1 + P^\mathrm{f} \preceq 0.
\end{align*}
Since $P_1\neq P^\mathrm{f}$ the Schur complement must necessarily have some negative eigenvalue. Consequently, no choice $-\Lambda \prec \bar P_\mathrm{s}$ can yield a positive semidefinite rotated stage cost matrix $H_\Lambda$. In turn, this entails that conditions~\eqref{eq:nl_conditions_lq} in their least restrictive form become
\begin{align*}
	P^\mathrm{f} \succ \bar P_\mathrm{s}, && H_\Lambda \succeq 0.
\end{align*}
which are only marginally different from the conditions we derive in this paper, which require the slightly stronger condition $H_\Lambda \succ 0.$ Note that this condition is required to make sure that the DARE does have a solution. In other words, this paper shows that the conditions in the companion paper~\cite{Gruene2024} for general nonlinear problems are almost tight in the nonsingular linear quadratic case. ``Almost'' here refers to the fact that there is no general statement for the case $P^\mathrm{f} \succeq \bar P_\mathrm{s}$.

As a side remark, we also observe that, using the same arguments as for $\bar P_\mathrm{s}$, one also obtains that no choice $-\Lambda \succ P_\mathrm{s}$ can yield a positive semidefinite rotated stage cost matrix $H_\Lambda$.

\subsection{Relation to First-Order Cost Corrections}
In this paper we focused on the quadratic term in the terminal cost, while in~\cite{Zanon2018a,Faulwasser2018} the impact of the linear term in the quadratic cost was analyzed. In particular, it has been proven that a wrong gradient (in the context of this paper the correct gradient is $\nabla V^\mathrm{f}(0)=0$) impedes asymptotic stability, though under some technical assumptions one still obtains practical stability~\cite{Grune2013a}. Even if the system starts at the optimal steady state (in our case the origin), the closed-loop system immediately leaves it. Our result allows us to further comment on this aspect. Indeed, in case the antistabilizing solution $P_\mathrm{a}$ exists, selecting $P^\mathrm{f}=P_\mathrm{a}$ yields a destabilizing controller for all initial states $x_0\neq0$. Clearly, this is a limit case of little practical interest, but it illustrates the fact that, while first-order conditions distinguish to which steady-state the system is stabilized, if at all, second-order conditions on the terminal cost, instead, make the difference between stability and instability, with the special limit case discussed above.

\section{Simulation Results}
\label{sec:simulations}

\begin{figure}
\includegraphics[width=\linewidth]{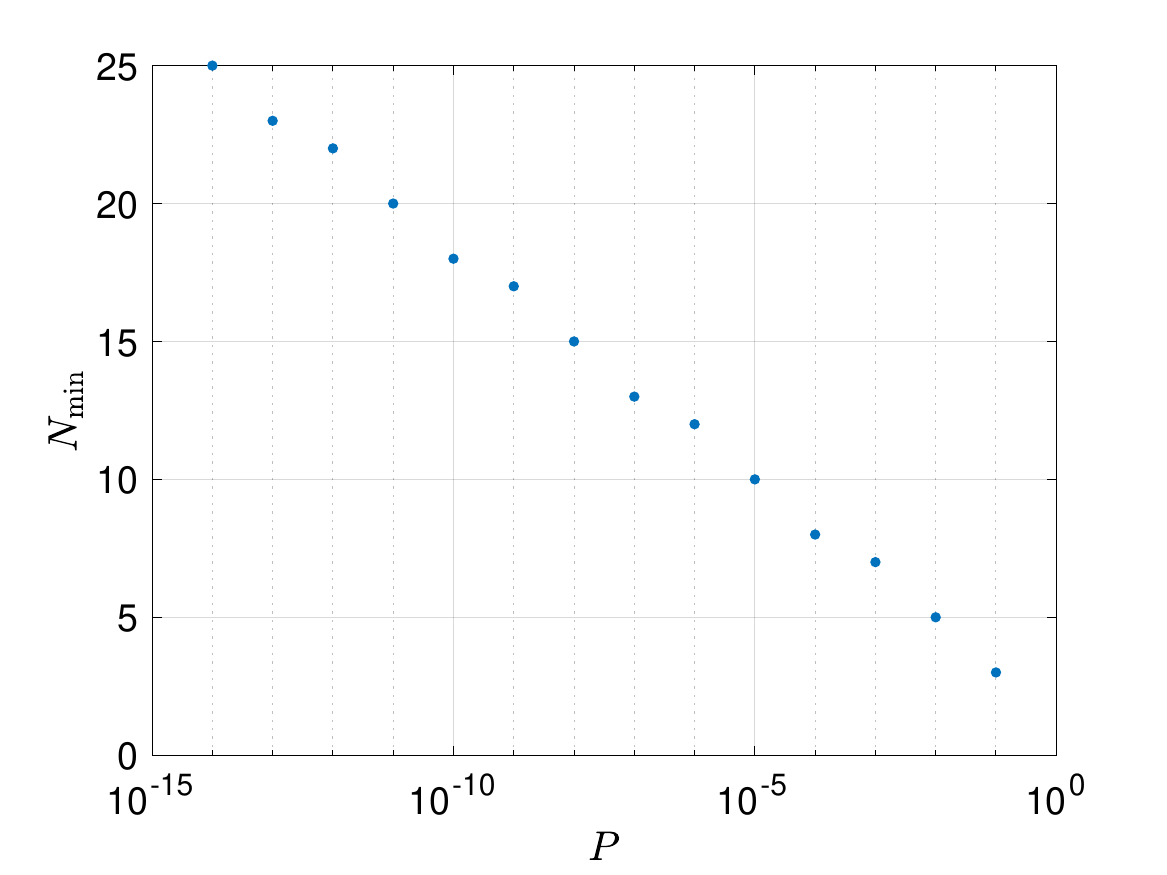}
\caption{Minimum stabilizing horizon length.}
\label{fig:ex1_horizon}
\end{figure}

\begin{figure}
	\includegraphics[width=\linewidth]{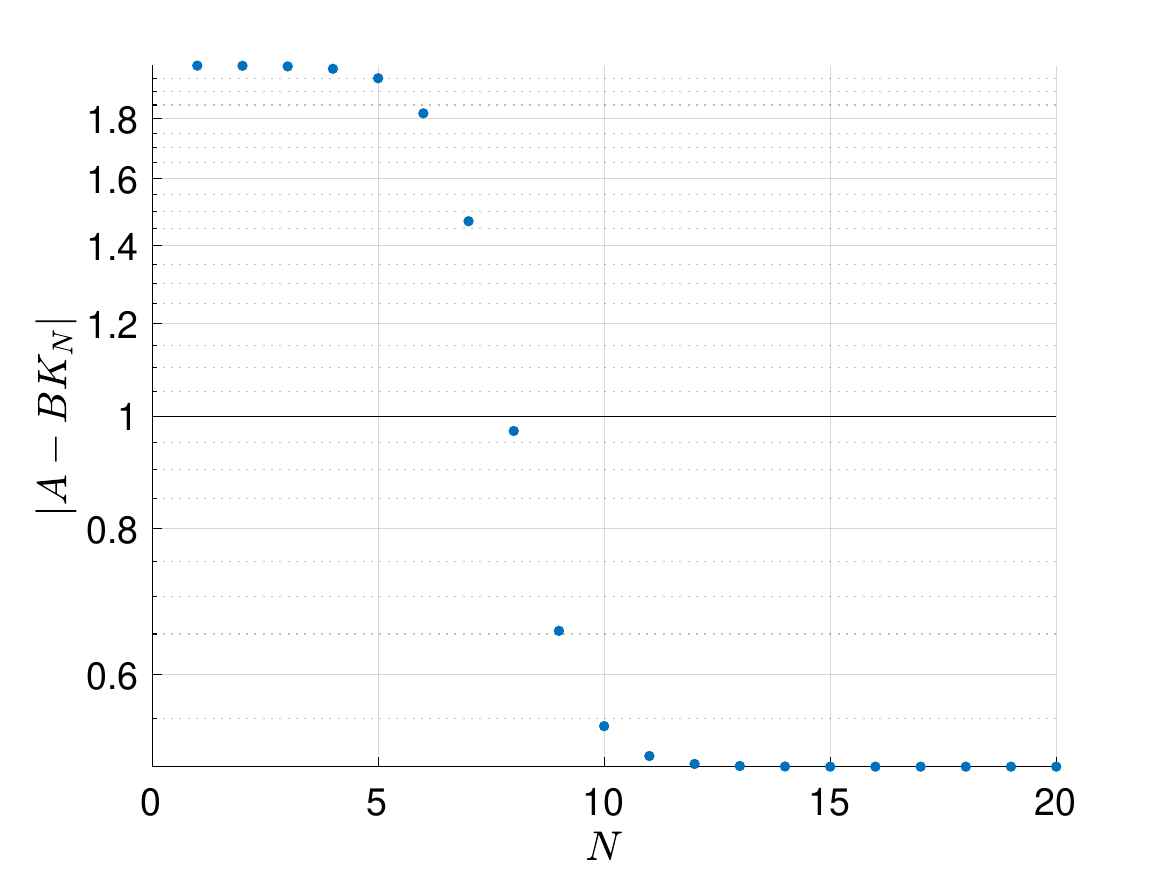}
	\caption{Closed-loop matrix eigenvalues with $P^\mathrm{f}=10^{-4}$.}
	\label{fig:ex1_ABK}
\end{figure}

Consider the case of Example~\ref{ex:simple_lqr2}. We recall that $P_\mathrm{a}=0$ and we display in Figure~\ref{fig:ex1_horizon} the minimum horizon length yielding asymptotic stability, which we computed for several terminal costs by solving the Riccati equation and evaluating the eigenvalues of $A-BK_N$. Furthermore, we display in Figure~\ref{fig:ex1_ABK} the absolute value of the eigenvalue of the closed-loop systme matrix for different $N$ and $P^\mathrm{f}=10^{-4}$. One can see that stability is obtained for $N\geq8.$

\begin{figure}
	\includegraphics[width=\linewidth]{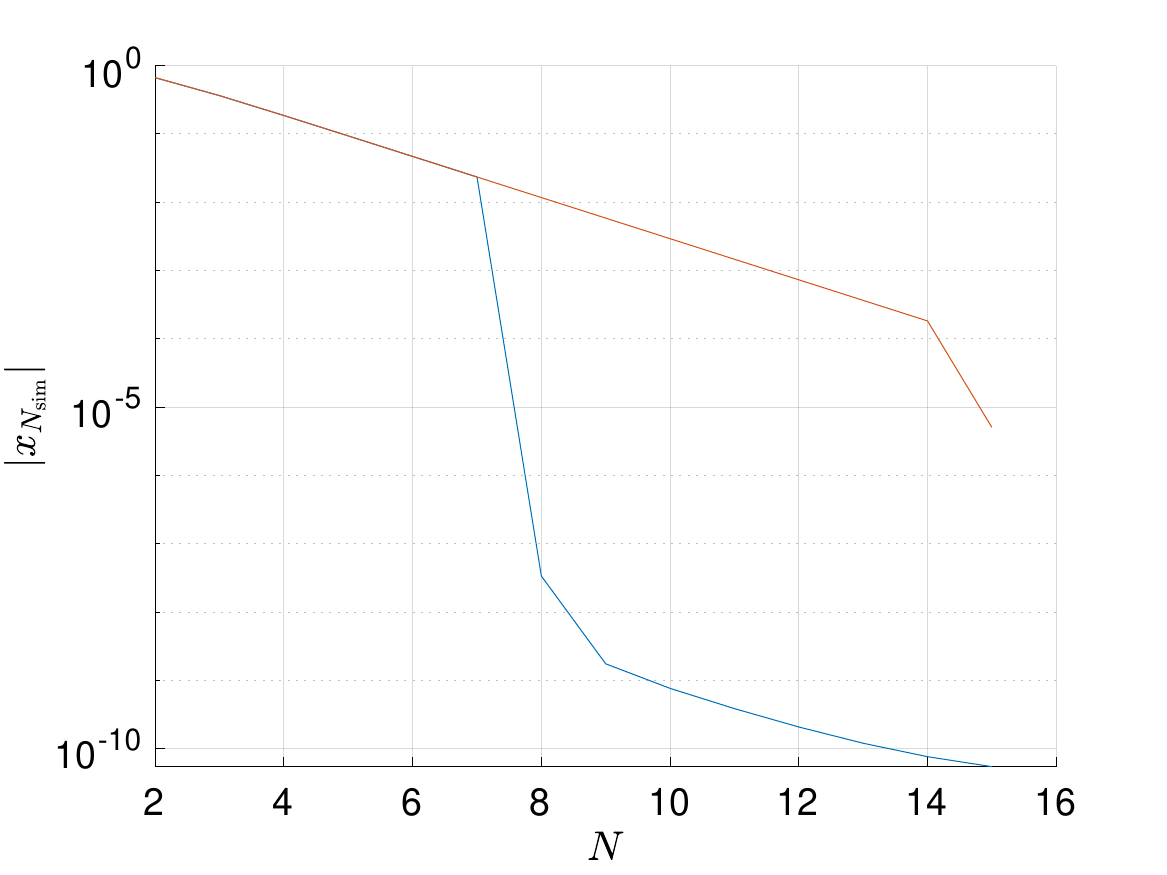}
	\caption{State at time $k=500$ obtained with $P^\mathrm{f}=0$ (red) and $P^\mathrm{f}=1e-4$ (blue).}
	\label{fig:ex1_xss}
\end{figure}

Let us introduce the constraint $|x|\leq1$. Since this constraint ensures that the state remains bounded, also the storage function remains bounded and the results of~\cite{Grune2013a} guarantee that the system converges to a neighborhood of the origin whose size decreases with increased prediction horizons. The results of~\cite{Faulwasser2018} guarantee convergence to the origin for a sufficiently long prediction horizon, under the condition that the unconstrained control law associated with the local linear-quadratic approximation of the system is stabilizing. For this example, unless one adds a terminal cost the unconstrained control law is not stabilizing. Indeed, we observe that in this setting, once the prediction horizon is long enough such that the unconstrained control law becomes stabilizing, then the system is indeed stabilized to the origin. This is shown in Figure~\ref{fig:ex1_xss}, where, after starting from initial state $\hat x_0=1$, 
the final state after $N_\mathrm{sim}=500$ time steps is never exactly $0$ is due to, on the one hand the finite simulation horizon, and, on the other hand, the fact that the MPC QP is solved to finite precision. Nevertheless, a clear difference between the case $P^\mathrm{f}=10^{-4}$ and the case $P^\mathrm{f}=0$ is seen for $N>8$, which is the horizon length starting from which the LQR is stabilizing also in the absence of constraints, as shown in Figure~\ref{fig:ex1_ABK}.

\section{Conclusions}
\label{sec:conclusions}

We have discussed the role of suitably chosen terminal costs in yielding exponential stability in the linear-quadratic case by establishing a strong connection to the symmetric solutions of the associated (reverse) discrete-time algebraic Riccati equation. We have connected our results to those obtained for the full nonlinear case and we have provided simple examples to illustrate our results.

\bibliographystyle{IEEEtran}
\bibliography{bibliography,bibliography_lars}

\begin{thebibliography}{10}
\providecommand{\url}[1]{#1}
\csname url@samestyle\endcsname
\providecommand{\newblock}{\relax}
\providecommand{\bibinfo}[2]{#2}
\providecommand{\BIBentrySTDinterwordspacing}{\spaceskip=0pt\relax}
\providecommand{\BIBentryALTinterwordstretchfactor}{4}
\providecommand{\BIBentryALTinterwordspacing}{\spaceskip=\fontdimen2\font plus
\BIBentryALTinterwordstretchfactor\fontdimen3\font minus
  \fontdimen4\font\relax}
\providecommand{\BIBforeignlanguage}[2]{{%
\expandafter\ifx\csname l@#1\endcsname\relax
\typeout{** WARNING: IEEEtran.bst: No hyphenation pattern has been}%
\typeout{** loaded for the language `#1'. Using the pattern for}%
\typeout{** the default language instead.}%
\else
\language=\csname l@#1\endcsname
\fi
#2}}
\providecommand{\BIBdecl}{\relax}
\BIBdecl

\bibitem{Rawlings2017}
J.~B. Rawlings, D.~Q. Mayne, and M.~Diehl, \emph{{M}odel {P}redictive
  {C}ontrol: {T}heory, {C}omputation, and {D}esign}, 2nd~ed.\hskip 1em plus
  0.5em minus 0.4em\relax Nob Hill Publishing, 2017.

\bibitem{Gruene2017}
L.~Gr{\"u}ne and J.~Pannek, \emph{{N}onlinear {M}odel {P}redictive {C}ontrol},
  2nd~ed., ser. Communications and Control Engineering.\hskip 1em plus 0.5em
  minus 0.4em\relax Springer International Publishing, 2017.

\bibitem{DiAR10}
M.~Diehl, R.~Amrit, and J.~B. Rawlings, ``A {L}yapunov function for economic
  optimizing model predictive control,'' \emph{IEEE Trans. Autom. Control},
  vol.~56, pp. 703--707, 2011.

\bibitem{Amrit2011a}
R.~Amrit, J.~Rawlings, and D.~Angeli, ``{E}conomic optimization using model
  predictive control with a terminal cost,'' \emph{Annual Reviews in Control},
  vol.~35, pp. 178--186, 2011.

\bibitem{Zanon2018a}
M.~Zanon and T.~Faulwasser, ``{E}conomic {MPC} without terminal constraints:
  {G}radient-correcting end penalties enforce asymptotic stability,''
  \emph{Journal of Process Control}, vol.~63, pp. 1 -- 14, 2018.

\bibitem{Willems1972a}
J.~Willems, ``\BIBforeignlanguage{English}{{D}issipative {D}ynamical {S}ystems
  {P}art {I}: {G}eneral {T}heory},'' \emph{\BIBforeignlanguage{English}{Archive
  for Rational Mechanics and Analysis}}, vol.~45, no.~5, pp. 321--351, 1972.

\bibitem{Willems1972b}
------, ``\BIBforeignlanguage{English}{{D}issipative {D}ynamical {S}ystems
  {P}art {II}: {L}inear {S}ystems with {Q}uadratic {S}upply {R}ates},''
  \emph{\BIBforeignlanguage{English}{Archive for Rational Mechanics and
  Analysis}}, vol.~45, no.~5, pp. 352--393, 1972.

\bibitem{Willems1971}
J.~C. {W}illems, ``{L}east {S}quares {S}tationary {O}ptimal {C}ontrol and the
  {A}lgebraic {R}iccati {E}quation,'' \emph{IEEE Transactions on Automatic
  Control}, vol. AC-16, no.~6, pp. 621--634, 1971.

\bibitem{DGSW14}
T.~Damm, L.~Gr\"une, M.~Stieler, and K.~Worthmann, ``An exponential turnpike
  theorem for dissipative discrete time optimal control problems,'' \emph{SIAM
  J. Control Optim.}, vol.~52, no.~3, pp. 1935--1957, 2014.

\bibitem{Gruene2018}
L.~Gr\"une and R.~Guglielmi, ``Turnpike properties and strict dissipativity for
  discrete time linear quadratic optimal control problems,'' \emph{SIAM Journal
  on Control and Optimization}, vol.~56, no.~2, pp. 1282--1302, 2018.

\bibitem{Gruene2024}
\BIBentryALTinterwordspacing
L.~Grüne and M.~Zanon, ``Stabilization of strictly pre-dissipative nonlinear
  receding horizon control by terminal costs,'' 2024. [Online]. Available:
  \url{https://arxiv.org/abs/2412.13538}
\BIBentrySTDinterwordspacing

\bibitem{Diehl2011}
M.~Diehl, R.~Amrit, and J.~Rawlings, ``{A} {L}yapunov {F}unction for {E}conomic
  {O}ptimizing {M}odel {P}redictive {C}ontrol,'' \emph{IEEE Trans. of Automatic
  Control}, vol.~56, no.~3, pp. 703--707, March 2011.

\bibitem{Angeli2012a}
D.~Angeli, R.~Amrit, and J.~Rawlings, ``{O}n {A}verage {P}erformance and
  {S}tability of {E}conomic {M}odel {P}redictive {C}ontrol,'' \emph{IEEE
  Transactions on Automatic Control}, vol.~57, pp. 1615 -- 1626, 2012.

\bibitem{Faulwasser2018a}
T.~Faulwasser, L.~Gr\"une, and M.~M{\"u}ller, ``Economic nonlinear model
  predictive control: Stability, optimality and performance,''
  \emph{Foundations and Trends in Systems and Control}, vol.~5, no.~1, pp.
  1--98, 2018.

\bibitem{Muller2015}
M.~A. M{\"u}ller, D.~Angeli, and F.~Allg{\"o}wer, ``On necessity and robustness
  of dissipativity in economic model predictive control,'' \emph{IEEE
  Transactions on Automatic Control}, vol.~60, no.~6, pp. 1671--1676, 2015.

\bibitem{Muller2013a}
M.~M\"uller, D.~Angeli, and F.~Allg\"ower, ``{O}n convergence of averagely
  constrained economic {MPC} and necessity of dissipativity for optimal
  steady-state operation,'' in \emph{Proceedings of the American Control
  Conference}, 2013.

\bibitem{Faulwasser2018}
T.~Faulwasser and M.~Zanon, ``{A}symptotic {S}tability of {E}conomic {NMPC}:
  {T}he {I}mportance of {A}djoints,'' in \emph{Proceedings of the IFAC
  Nonlinear Model Predictive Control Conference}, 2018.

\bibitem{Grune2013a}
L.~Gr{\"u}ne, ``{E}conomic receding horizon control without terminal
  constraints,'' \emph{Automatica}, vol.~49, pp. 725--734, 2013.

\bibitem{Muller2016}
M.~A. M\"uller and L.~Gr\"une, ``Economic model predictive control without
  terminal constraints for optimal periodic behavior,'' \emph{Automatica},
  vol.~70, pp. 128 -- 139, 2016.

\bibitem{Molinari1975a}
B.~P. Molinari, ``{T}he {S}tabilizing {S}olution of the {D}iscrete {A}lgebraic
  {R}iccati {E}quation,'' \emph{IEEE Transactions on Automatic Control},
  vol.~20, pp. 396--399, 1975.

\bibitem{Lancaster1986}
P.~Lancaster, A.~C.~M. Ran, and L.~Rodman, ``{H}ermitian solutions of the
  discrete algebraic {R}iccati equation,'' \emph{International Journal of
  Control}, vol.~44, pp. 777--802, 1986.

\bibitem{Ran1993a}
A.~C.~M. Ran and H.~L. Trentelman, ``{L}inear {Q}uadratic {P}roblems with
  {I}ndefinite {C}ost for {D}iscrete {T}ime {S}ystems,'' \emph{SIAM Journal on
  Matrix Analysis and Applications}, vol.~14, pp. 776--797, 1993.

\bibitem{Gruene2014}
L.~Gr{\"u}ne and M.~Stieler, ``A {L}yapunov function for economic {MPC} without
  terminal conditions,'' in \emph{Proc. of the 53rd IEEE Conference on Decision
  and Control}, December 2014, pp. 2740--2745.

\bibitem{Zanon2014d}
M.~Zanon, S.~Gros, and M.~Diehl, ``{I}ndefinite {L}inear {MPC} and
  {A}pproximated {E}conomic {MPC} for {N}onlinear {S}ystems,'' \emph{Journal of
  Process Control}, vol.~24, pp. 1273--1281, 2014.

\bibitem{Stoorvogel1998}
A.~A. Stoorvogel and A.~Saberi, ``The {D}iscrete {A}lgebraic {R}iccati
  {E}quation and {L}inear {M}atrix {I}nequality,'' \emph{Linear Algebra and its
  Applications}, vol. 274, no.~1, pp. 317--365, 1998.

\bibitem{Ionescu1996}
V.~Ionescu, ``Reverse discrete-time riccati equation and extended nehari's
  problem,'' \emph{Linear Algebra and its Applications}, vol. 236, pp. 59--94,
  1996.

\bibitem{Zanon2016b}
M.~Zanon, S.~Gros, and M.~Diehl, ``{A} {T}racking {MPC} {F}ormulation that is
  {L}ocally {E}quivalent to {E}conomic {MPC},'' \emph{Journal of Process
  Control}, vol.~45, pp. 30 -- 42, 2016.

\bibitem{Anderson1990}
B.~D.~O. Anderson and J.~B. Moore, \emph{{Optimal Control - Linear Quadratic
  Methods}}.\hskip 1em plus 0.5em minus 0.4em\relax Dover, 1990.

\end{thebibliography}


\end{document}